%% file: Memoire.tex
\documentclass[11pt,a4paper,twoside,openright]{book}

\usepackage{multind}  
\makeindex{symboles}		

\usepackage{setspace}

\usepackage[utf8]{inputenc}
\usepackage[T1]{fontenc}
\usepackage[english]{babel}

\usepackage{amsmath,amssymb,dsfont}
\usepackage[mathscr]{euscript}
\usepackage{bm}
\usepackage{amsthm}  
\usepackage{wasysym}

\usepackage{wrapfig} 
\usepackage[textwidth=14.5cm,textheight=23cm]{geometry}
\usepackage{pdfsync}                  
\usepackage{boites}
\usepackage{chngpage}
\usepackage{twcal}
\usepackage[T1]{fontenc}

\usepackage{pstricks}

\usepackage[bookmarks=true,hyperfootnotes=false]{hyperref}
\hypersetup{
			colorlinks=true,
			linkcolor=oxfordblue,
			anchorcolor=black,
			citecolor=black,
			urlcolor=persianindigo,
			pdftitle={Algorithms, unaffected by the Schwarz paradox},
			pdfauthor={Paolo Roselli},
			pdfkeywords={Schwarz paradox, Inscribed polyhedra, Tangent planes, Surface area, Geometric Algebra, Clifford Algebra}
}

\definecolor{ccita}{rgb}{1, 0.84, 0.6}
\definecolor{dgreen}{rgb}{0,0.4,0} 
\definecolor{gold}{rgb}{0.85,.66,0}
\definecolor{armygreen}{rgb}{0.29, 0.33, 0.13} 
\definecolor{airforceblue}{rgb}{0.36, 0.54, 0.66}
\definecolor{arsenic}{rgb}{0.23, 0.27, 0.29}
\definecolor{ao}{rgb}{0.0, 0.0, 1.0}
\definecolor{oxfordblue}{rgb}{0.0, 0.13, 0.28}
\definecolor{persianindigo}{rgb}{0.2, 0.07, 0.48}
\definecolor{auburn}{rgb}{0.43, 0.21, 0.1}

\theoremstyle{plain}

\newtheorem{definition}{D\'efinition}[chapter]
\newtheorem*{definition*}{D\'efinition}

\newtheorem*{quest*}{Question}

\newtheorem*{nota*}{Notation}
\newtheorem{exam}[definition]{Example}
\newtheorem*{exam*}{Example}

\newtheorem*{exer*}{Exercice}
\newtheorem{rem}[definition]{Remark}
\newtheorem*{rem*}{Remark}

\newtheorem*{cita*}{Citation}
\newtheorem{proposition}[definition]{Proposition}
\newtheorem*{proposition*}{Proposition}
\newtheorem{theorem}[definition]{Theorem}
\newtheorem*{theorem*}{Theorem}

\newtheorem*{coro*}{Corollaire}
\newtheorem{lem}[definition]{Lemme}
\newtheorem*{lem*}{Lemme}

\newtheorem*{conj*}{Conjecture}

\newsavebox{\fmbox}

\newenvironment{defi*}{%
	\begin{definition*}
	\mbox{}
	\begin{breakbox}
	\begin{adjustwidth}{1cm}{0.5cm}
	}{
	\end{adjustwidth}
	\end{breakbox}
    \end{definition*}
}

\newenvironment{prop*}{%
	\begin{proposition*}
	\mbox{}
	\begin{breakbox}
	\begin{adjustwidth}{1cm}{0.5cm}
	}{
	\end{adjustwidth}
	\end{breakbox}
    \end{proposition*}
}

\newenvironment{theo*}{%
	\begin{theoreme*}
	\mbox{}
	\begin{breakbox}
	\begin{adjustwidth}{0.7cm}{1cm}
	}{
	\end{adjustwidth}
	\end{breakbox}
    \end{theoreme*}
}

\newenvironment{lemme*}{%
	\begin{lem*}
	\mbox{}
	\begin{breakbox}
	\begin{adjustwidth}{0.7cm}{1cm}
	}{
	\end{adjustwidth}
	\end{breakbox}
    \end{lem*}
}

\def\1{\mathds{ 1}}

\title{\bf Algorithms, unaffected by the Schwarz paradox, approximating tangent planes and area of smooth surfaces via inscribed triangular polyhedra\\ \ } 
\author{{\Large\bf Paolo Roselli} \\ \\ Dipartimento di Matematica \\ Universit\`a di Roma ``Tor Vergata'' \\ Via della Ricerca Scientifica \\ 00133 - Rome \\ ITALY \\ \textsc{roselli@mat.uniroma2.it} \\ \\
Institut de Recherche en Math\'ematique et Physique \\
Universit\'e Catholique de Louvain \\ Chemin du Cyclotron 2 \\ 1348 - Louvain la Neuve \\ BELGIUM \\ \textsc{paolo.roselli@uclouvain.be}
\\ \ \\ \ \\ \ \\ \ \\ \ \\ \ \\ 
M\'emoire pr\'esent\'e au concours annuel de 2014 \\ de la Classe des Sciences 
de l'Acad\'emie Royale de Belgique, \\ Groupe I - Math\'ematiques: d, portant sur \\ une contribution aux applications de l'alg\'ebre de Clifford \`a l'analyse.}

\date{}

\begin{document}

\maketitle

\newpage
\thispagestyle{empty}
\mbox{}
\setcounter{page}{0}


\tableofcontents

\input{01-Introduction.tex}

\input{02-GeometricAlgebra.tex}

\input{03-Geometry.tex}

\input{04-SmoothCurves.tex}

\input{05-SmoothSurfaces.tex}

\input{06-MainResults.tex}

\input{07-LocalSchwarzParadox.tex}

\input{09-Bibliography.tex}
\cleardoublepage
\phantomsection
\printindex{symboles}{Symbols}

\cleardoublepage
\phantomsection
\printindex{termes}{Index}

\end{document}

%% file: 01-Introduction.tex
\chapter{Introduction}

\section{Goals}

In this work we provide algorithms\footnote{See~(\ref{eq:AlgI}) in Theorem~\ref{thm:II} and~(\ref{eq:AlgII}) in Theorem~\ref{thm:III}.} approximating the bivector $\displaystyle \partial_{\ell_1} s_{(x)} \wedge \partial_{\ell_2} s_{(x)}$ and the integral $\displaystyle \int_P \big|\partial_{\ell_1} s(x) \wedge \partial_{\ell_2} s(x) \big|dx$ of a smooth map $s:\Omega \to \mathbb{E}_n$ (that we loosely call `surface'), where

\begin{itemize}
	\item {\color{dgreen} $\mathbf{\mathbb{E}_n}$}\index{symboles}{$\mathbb{E}_n$} is an $n$-dimensional Euclidean space;
	\item $\Omega$ is an open subset of the Euclidean plane $\mathbb{E}_2$;
	\item $P\subset\Omega$ is a compact polygon;
	\item for every $v\in \mathbb{E}_2$, $\displaystyle \partial_v s_{(x)}= \lim_{\epsilon \to 0}\frac{1}{\epsilon} \big[s(x+\epsilon v)-s(x)\big]$;
	\item $\{\ell_1, \ \ell_2\}$ is an orthonormal basis in $\mathbb{E}_2$;
	\item $\wedge$ is the outer product in the Euclidean Clifford algebra $\mathbb{G}_n$ associated\footnote{See Sections~\ref{sec:Eucl struct} and~\ref{sec:GA and En}.} to $\mathbb{E}_n$.
\end{itemize}

In particular, if $\displaystyle \partial_{\ell_1} s_{(x)} \wedge \partial_{\ell_2} s_{(x)}\ne 0$, then the bivector $\displaystyle \partial_{\ell_1} s_{(x)} \wedge \partial_{\ell_2} s_{(x)}$ can represent\footnote{See Section~\ref{sec:point lines}} the direction of the tangent plane to the surface $s$ at point $s(x)$ (or the normal vector, if $s:\Omega \to \mathbb{E}_3$ and if we consider\footnote{See Section~\ref{sec:cross product}.} the cross product $\partial_{\ell_1} s_{(x)} \bm{\times} \partial_{\ell_2} s_{(x)}$).

Our algorithms use informations from triangles in $\mathbb{E}_n$ inscribed\footnote{This means that the vertices of the triangles are images $s(x)$ of vertices of some nondegenerate triangles in $\Omega$ (see also Section~\ref{sec:surfaces}).} in the surface~$s$. Thus, Algorithm~(\ref{eq:AlgI}) allows to recover the tangent plane direction from every sequence of inscribed triangles converging to the point $s(x)$; this result is obtained approximating\footnote{See~(\ref{eq:Alg0}) in Theorem~\ref{thm:I}.} Jacobian determinants of smooth transformations $f:\Omega \to \mathbb{E}_2$ at points $x\in \Omega$ through nondegenerate triangles converging to point $x$. Algorithm~(\ref{eq:AlgI}) can also estimate the norm of $\displaystyle \partial_{\ell_1} s_{(x)} \wedge \partial_{\ell_2} s_{(x)}$, and thus, when $s$ is globally injective, Algorithm~(\ref{eq:AlgII}) can approximate the area of portions of the surface $s$ from every sequence of inscribed triangular\footnote{This means that all faces of the polyhedron are inscribed triangles.} polyhedra uniformly convergent to that portion\footnote{See Remark~\ref{rem:triangulations}.}.

In particular, we apply Algorithm~(\ref{eq:AlgI}) to the triangulation of a circular cylinder of the famous Schwarz\footnote{Hermann Amandus Schwarz (1843-1921).} area paradox\footnote{See Chapter~\ref{cha:local Schwarz}.}, showing that the approximating inscribed balanced mean bivectors\footnote{See Section~\ref{sec:surfaces}.} do converge to the tangent bivectors without any restriction of the approximating triangular mesh. 

As a matter of fact, by using Algorithms~(\ref{eq:AlgI}) and~(\ref{eq:AlgII}) we can restore analogies\footnote{Compare, for instance, Proposition~\ref{prop: approxim dot c} and Theorem~\ref{thm:II}.} between the limit vector $\dot{c}_{(\chi)}$ of a smooth curve $c:I\to \mathbb{E}_n$ and the limit bivector $\displaystyle \partial_{\ell_1} s_{(x)} \wedge \partial_{\ell_2} s_{(x)}$ of a smooth surface $s:\Omega\to \mathbb{E}_n$; such analogies are lost, according to the Schwarz paradox, if we try to approximate tangents or surface area via the usual algorithms applied to arbitrary inscribed triangular polyhedra.

\section{Warnings}\label{sec:warnings}

The aim of this work is to describe Algorithms~(\ref{eq:AlgI}) and~(\ref{eq:AlgII}) as simply as possible; thus, our intention here is not to provide the most general hypothesis under which such algorithms work; neither do we want to generalize them here to $k$-manifolds immersed in $n$-dimensional Euclidean spaces, or to Riemann manifolds, nor do we want to introduce a Stieltjes-like $k$-measure in $\mathbb{E}_n$ generalizing Theorem~\ref{thm:I}. Such generalizations will be examined in forthcoming works.

The main theorems are stated and proved using Geometric Algebra. However, the reader will be provided formulas to translate them into the lengthy Cartesian coordinate formalism.  

Finally, we apologize if some  calculations may appear tedious or pedantic to  readers well acquainted with Geometric Algebra, but this work is addressed to a broader audience.

\section{Notations I}

In this work we consider it important to distinguish the different types of mathematical objects in our formulas; therefore, we use the following conventions:

\begin{itemize}
	\item lower-case Greek letters stand for real numbers;
	\item lower-case Latin letters stand for vectors in some Euclidean space $\mathbb{E}_n$ (with the exceptions of letters $i,\ j,\ k,\ m,\ n$, representing integer indexes);
	\item capital Latin letters stand for bivectors or generic $k$-vectors;
	\item capital Greek letters stand for sets;
	\item capital bold Greek letters stand for $n$-uples or arrays of real numbers (with $n>1$). 
\end{itemize}

\section{Historical notes}

As two distinct points on a sufficiently smooth curve converge to the same point, the line passing through those two points assumes a well defined position. In particular, when such a local phenomenon is globally injective and uniform, we can approximate the length of the curve by the lengths of the line segments joining a finite number of consecutive points on the curve. The idea that a similar phenomenon may occur to triangles inscribed in a sufficiently smooth surface is probably what suggested to Serret\footnote{Joseph Alfred Serret (1819-1885).} (see \cite{Ser1879}) the following definition of area\footnote{Our translation: ``Let a portion of a curved surface be bounded by a contour C; we will call area of that surface the limit S to which converges the area of an inscribed polyhedral surface whose faces are triangles and which is bounded by a polygonal contour F having C as limit.''}:
\bigskip

\textit{\large Soit une portion de surface courbe termin\'ee par un contour C; nous nommerons aire de cette surface la
limite S vers laquelle tend l'aire d'une surface poly\'edrale inscrite form\'ee de faces triangulaires et termin\'ee par un contour polygonal F ayant pour limite le contour~C.}
\bigskip

However, on 20 December 1880, Schwarz wrote to Genocchi\footnote{Angelo Genocchi (1817-1889).} (see \cite{Cas1950}) observing that the area of a curved surface cannot be defined as Serret did. In subsequent letters to Genocchi, Schwarz showed that even the area of a surface as simple as a bounded part of a right circular cylinder cannot be recovered using Serret's definition.
Schwarz even provided examples of sequences of inscribed triangular polyhedra whose areas converge to any given number not less than the area of the cylinder (and even to infinity) as the polyhedra approach uniformly the cylinder\footnote{As a consequence, there also exist sequences of inscribed triangular polyhedra approaching the cylinder whose areas have no limit.}. Such phenomenon, that may occur to every curved surface (and even to polyhedra\footnote{See \cite{Fre1925}.}) is given the name of {\bf \color{dgreen} Schwarz paradox}
\index{termes}{Schwarz paradox} (or {\bf \color{dgreen} Schwarz phenomenon}).
\index{termes}{Schwarz phenomenon}

That famous paradox apparently destroyed the possibility of defining the area of a smooth surface by analogy with the length of a smooth curve. 
Besides, the local interpretation of the Schwarz phenomenon implies that as three noncollinear points on a smooth surface converge to the same point of the surface, the limit position of the plane passing through those three points is not well determined, and can differ from the tangent plane to the surface at the limiting point.
Also, Schwarz's counterexample shows that the limiting position of the secant plane can even be orthogonal to the actual tangent plane.
\bigskip

Two questions naturally arise:

\begin{itemize}
	\item what sequences of inscribed triangular polyhedra approaching a surface  have areas converging to the area of that surface~?
	\item are there algorithms able to recover the area of a surface from every sequence of inscribed triangular polyhedra approaching that surface~?
\end{itemize}

Schwarz showed that those questions are not trivial even for a cylinder.
\bigskip

Many different approaches were used to answer those questions. We cannot summarize such a long and prolific history here\footnote{Suggested readings are~\cite{GanPer2009}, \cite{Ces1956}, \cite{Rad1948} and~\cite{Sak1937}.}; we will just focus on some particular issues which only concern smooth curved surfaces.

\begin{itemize}
	\item 
	Apart from Peano\footnote{Giuseppe Peano (1858-1932).}, all authors\footnote{Another exception is William Henry Young (1863-1942), who used an approach similar to Peano's in~\cite{You1920}.} approximated the area of a curved surface using the areas of triangular polyhedra uniformly approaching the surface.
	\begin{itemize}
		\item 
	Most of those authors selected particular inscribed triangular polyhedra constraining the form or the position of the triangular faces with `ad hoc' conditions. 
		\item 
	Lebesgue\footnote{Henri Lebesgue (1875-1941).}, instead, freed himself from inscribed polyhedra and artificial geometric conditions; however, his definition of area\footnote{See~\cite{Leb1902}} is of no help in selecting a sequence of polyhedra whose areas converge to the area of the surface\footnote{We suggest to read the Jordan's criticism to Lebesgue's definition of area in \cite{Leb1926} at pages~163--164.}, nor does his definition of area correspond locally to a definition of tangent. 
	\item
Ge\"{o}cze\footnote{Zo\'ard Ge\"{o}cze (1873-1916).} conjectured\footnote{See \cite{Rad1948}.}, and Mulholland\footnote{H.~P.~Mulholland (we did not find any biographical data about him).} proved in~\cite{Mul1950}, that Lebesgue's area can also be obtained restricting Lebesgue's approach to inscribed polyhedra.
	\end{itemize}
\item
Peano freed himself from polyhedra\footnote{Nevertheless, his work has been a key inspiration to us.}, and used his \textit{Calcolo Geometrico} (based on Grassmann exterior algebra\footnote{See~\cite{Pea1887} and~\cite{Pea1888}.}) to define the area through integrals taken on the boundaries of portions of a surface. However, his definition\footnote{See~\cite{Pea1887} on page~164, or \cite{Pea1890} on page~55.} was vague about what portions of a surface may cut in order to approximate the area of the whole surface.
\end{itemize}

Our Algorithm~(\ref{eq:AlgII}) allows to consider inscribed triangular polyhedra without any kind of constraint, and uses a slightly modified notion of area. Besides, Algorithm~(\ref{eq:AlgII}) is just a global adaptation of the local Algorithm~(\ref{eq:AlgI}) that approximates tangent planes from every inscribed triangle approaching a point on the surface. Thus, Algorithms~(\ref{eq:AlgI}) and~(\ref{eq:AlgII}) restore many of the analogies between curves and surfaces.

%% file: 02-GeometricAlgebra.tex
\chapter{Basic notions of Euclidean Clifford algebras}

\section{Motivations}

Theorem~\ref{thm:I}, Theorem~\ref{thm:II} and Theorem~\ref{thm:III} are stated and proved using Euclidean Clifford algebra (i.e. {\bf \color{dgreen} Geometric Algebra}\index{termes}{Geometric Algebra}). 
Of course, they can be translated into the Cartesian coordinatewise language as well\footnote{We provide formulas to do it.}; however, we consider the coordinate-free language of Geometric Algebra to be richer and more suitable in order to algebraically represent geometric properties. Moreover, we discovered Algorithm~(\ref{eq:AlgI}) and~(\ref{eq:AlgII}) while exploring the Schwarz paradox via Geometric Algebra and not via Cartesian language.

\section{Formal Geometric Algebra}

Following is a brief formal description of Geometric Algebra $\mathbb{G}_n$. 
For more details and other approaches, see also~\cite{Art2006},~\cite{Cli1876},~\cite{Del1992},~\cite{DorFonMan2007},~\cite{HesSob1984},~\cite{Hes1999},~\cite{Lou2001},~\cite{Mac2002} or~\cite{Sny2012}. 
\medskip

Suppose we have an ordered alphabet of $n$ (distinct) letters {\bf \color{dgreen} $\mathcal{A}_n$}$=\{\ell_1,\dots , \ell_n\}$\index{symboles}{$\mathcal{A}_n$}.
A {\bf \color{dgreen} word}\index{termes}{Word} from this alphabet is a juxtaposition of letters taken from $\mathcal{A}_n$. A word with no letters is considered a word as well, it is named {\bf \color{dgreen} empty word}\index{termes}{Empty word} (or {\bf \color{dgreen}unit}),\index{termes}{Unit} and it is given the reserved\footnote{A symbol is called `reserved' if it can never be a letter of any alphabet.} symbol~{\bf \color{dgreen}$\1$})\index{symboles}{$\1$}.

The set of formal finite real linear combinations of words\footnote{We will write real coefficients on left of words.} from $\mathcal{A}_n$ forms a real algebra {\bf \color{dgreen} $\mathbb{G}_n$}\index{symboles}{$\mathbb{G}_n$} if we consider juxtaposition of words as an associative and distributive product among words\footnote{The real coefficients are multiplied among themselves in $\mathbb{R}$.}. Thus, $\1$ is the unit for that product.

Also the empty real linear combination of words is considered an element of such algebra, and it is given the symbol~{\bf \color{dgreen}$\mathbb{ O}$}.\index{symboles}{$\mathbb{ O}$}
The following axioms hold in $\mathbb{G}_n$:

\begin{center}
$
\displaystyle 
\hfil
\ell_j \ne \1\ ,
\hfil
\ell_j \ne \mathbb{O}\ ,
\hfil
\mathbb{O} \ne \1\ ,
\hfil
0 W = \mathbb{O}\ ,
\hfil
1 W = W\ ,
\hfil
$
\end{center}

where $j=1,\dots, n$ and $W$ is a word from the alphabet $\mathcal{A}_n$; moreover,

\begin{equation}\label{eq: product axiom}
\ell_i \ell_j =
\left\{
\begin{array}{ll}
	-\ell_j \ell_i & \textrm{if } i\ne j \ , \\
	\1 & \textrm{if } i=j \ , 
\end{array}
\right.
\end{equation}

where $-\ell_j \ell_i$ abbreviates $(-1)\ell_j \ell_i$.
The complete ordered word $\ell_1 \ell_2 \cdots \ell_n$ is called {\bf \color{dgreen} pseudo-unit}\index{termes}{Pseudo-unit} and is given the reserved symbol~{\bf \color{dgreen}$\mathbb{I}_n$}\index{symboles}{$\mathbb{I}_n$}.
$\mathbb{G}_n$ is then uniquely determined\footnote{Unique up to algebra isomorphisms between real associative algebras with unit.} if we add the final axioms

\begin{center}
$\displaystyle
\hfil
\mathbb{I}_n \ne \mathbb{O}\ ,
\hfil
\mathbb{I}_n \ne \1\ ,
\hfil
\mathbb{I}_n \ne -\1 \ .
\hfil
$
\end{center}

Axiom~(\ref{eq: product axiom}) allows to reduce every nonempty word from the alphabet $\mathcal{A}_n$ to a unique minimal\footnote{That is, without repeated letters.} ordered word (with sign)

\[
\pm
\ell_{i_1} \ell_{i_2} \cdots \ell_{i_k}\ ,
\] 

where $i_1< i_2< \cdots <i_k$. The number $k$ is called {\bf \color{dgreen} grade}\index{termes}{Grade of a word} of the word, and the sign is called {\bf \color{dgreen} orientation}\index{termes}{Orientation of a word} of the word (with respect to the ordered alphabet $\mathcal{A}_n$). Such reductions make $\mathbb{G}_n$ a graded algebra

\[
\mathbb{G}_n
=
\bigoplus_{k=0}^n \mathbb{G}_{n \choose k}\ ,
\]

where~{\bf \color{dgreen}$\displaystyle \bm{\mathbb{G}_{n \choose k}}$}\index{symboles}{$\mathbb{G}_{n \choose k}$}
 is the linear subspace of (finite) real  combinations of words of grade $k$ (notice that $\displaystyle \mathbb{G}_{n \choose k}$ is not a subalgebra). Each $\displaystyle \mathbb{G}_{n \choose k}$ has (real) dimension ${n \choose k} = \frac{n!}{k!(n-k)!}$, and $\mathbb{G}_n$ has dimension $2^n$.

We can unambiguously identify the algebra of real numbers $\mathbb{R}$ with $\displaystyle \mathbb{G}_{n \choose 0}$, the real number $1$ with unit $\1$, and $0\in \mathbb{R}$ with the empty linear combination $\mathbb{O}$.

\section{Euclidean structure of $\mathbb{G}_n$}\label{sec:Eucl struct}

Geometric Algebra $\mathbb{G}_n$ is also called {\bf \color{dgreen} Euclidean Clifford algebra}\index{termes}{Euclidean Clifford algebra}, because it possesses a Euclidean structure strictly tied with its algebraic product\footnote{And because it was introduced by William Kingdon Clifford (1845-1879); see~\cite{Cli1876}.}. As a matter of fact, the symmetric part of the product among elements $x,y \in \mathbb{G}_{n \choose 1}$
  
\[
\frac{1}{2}
(xy+yx)
\]	

\noindent
is always a real number and, as a function of $x$ and $y$, it is a symmetric, positive definite bilinear form in $\mathbb{G}_{n \choose 1}$ (that we denote with the symbol~{\bf \color{dgreen} $\bm{x\cdot y}$}).\index{symboles}{$x\cdot y$}

The $n$ letters of the ordered alphabet $\mathcal{A}_n$ (generating $\mathbb{G}_n$) form an ordered orthonormal basis in $\mathbb{G}_{n \choose 1}$ with respect to the scalar product $x \cdot y$, indeed

\[
\ell_i \cdot \ell_j =
\frac{1}{2}
(\ell_i \ell_j + \ell_j \ell_i )
=
\left\{
\begin{array}{ll}
	0 & \textrm{if } i\ne j \ , \\
	1 & \textrm{if } i=j \ . 
\end{array}
\right.
\]

It is also important to note that the antisymmetric part of the product between $x,y \in \mathbb{G}_{n \choose 1}$

\[
\frac{1}{2}
(xy-yx)
\]

\noindent
is always an element of $\mathbb{G}_{n \choose 2}$, it is given the symbol~{\bf \color{dgreen} $\bm{x\wedge y}$}\index{symboles}{$x\wedge y$}, and is called {\bf \color{dgreen} outer product}\footnote{Or, simply, exterior product.}\index{termes}{Outer product in $\mathbb{G}_n$} because it acts in the graded algebra $\mathbb{G}_n$ as the associative antisymmetric Grassmann exterior product acts on the graded Grassmann algebra $\displaystyle \bigoplus_{k=0}^n \left[\bigwedge^{k} \mathbb{G}_{n \choose 1}\right]$. So, the product of $x,y \in \mathbb{G}_{n \choose 1}$ can be decomposed as 

\[
xy = (x\cdot y) \ + \ (x\wedge y)\ ,
\]

and 

\[
\ell_i \ell_j =
(\ell_i \cdot \ell_j) 
+
(\ell_i \wedge \ell_j) 
=
\left\{
\begin{array}{ll}
\ell_i \wedge \ell_j = - \ell_j \wedge \ell_i	 & \textrm{if } i\ne j \ , \\
1 & \textrm{if } i=j \ . 
\end{array}
\right.
\]

\section{The Geometric Algebra associated to an oriented Euclidean space}\label{sec:GA and En}

Here we notice that the construction of $\mathbb{G}_n$ can proceed the other way around as well: given a $n$-dimensional Euclidean space $\mathbb{E}_n$, 
there exists a unique\footnote{Up to isometries and orientation.} Geometric Algebra $\mathbb{G}_n$ such that its Euclidean subspace $\mathbb{G}_{n \choose 1}$ is isometric to $\mathbb{E}_n$. Indeed, it suffices to choose an ordered orthonormal basis $\{e_1,\ e_2,\ \dots ,\ e_n\}\subset \mathbb{E}_n$ as the ordered alphabet generating $\mathbb{G}_n$. In this sense we speak of Geometric Algebra associated to the oriented\footnote{The orientation being determined by the order of the orthonormal basis.} Euclidean space~$\mathbb{E}_n$.

A fundamental improvement of Geometric Algebra  $\displaystyle \mathbb{G}_n
=\bigoplus_{k=0}^n \mathbb{G}_{n \choose k}$ over Grassmann algebra $\displaystyle \bigoplus_{k=0}^n \left[\bigwedge^{k} \mathbb{E}_n\right]$ is that $\mathbb{G}_n$ has a well defined Euclidean structure such that

\begin{itemize}
	\item each subspace $\mathbb{G}_{n \choose k}$ is orthogonal in $\mathbb{G}_n$ to every other $\mathbb{G}_{n \choose j}$, with $k\ne j$;
	\item each subspace $\mathbb{G}_{n \choose k}$ has a Euclidean structure, i.e. a symmetric, positive-definite bilinear form (that we will continue to indicate with the dot $\bm{\cdot}$) uniquely determined by the scalar product in $\mathbb{G}_{n \choose 1}$ (usually identified with $\mathbb{E}_n$).
\end{itemize}

For instance, for each $a,b,c,d \in \mathbb{G}_{3 \choose 1}\equiv \mathbb{E}_3$

\[
(a\cdot c) (b\cdot d)- (a\cdot d)(b\cdot c)
\]

is the scalar product $(a\wedge b) \cdot (c\wedge d)$ in $\mathbb{G}_3$ restricted to the subspace $\mathbb{G}_{3 \choose 2}$.

Notice that the one-dimensional subspace $\mathbb{G}_{3 \choose 0}$ (that we identified with $\mathbb{R}$) has the usual product between real numbers as the restriction of the scalar product in~$\mathbb{G}_3$

\begin{center}
$
\big(\alpha\1\big)
\cdot
\big(\beta\1\big)
=
\alpha\beta
=
\big(\alpha\1\big)
\big(\beta\1\big)
$,
\end{center}

while

\begin{center}
$
\big(\alpha\mathbb{I}_n\big)
\cdot
\big(\beta\mathbb{I}_n\big)
=
\alpha\beta
$, and 
$
\big(\alpha\mathbb{I}_n\big)
\big(\beta\mathbb{I}_n\big)
=
(-1)^{\frac{n(n-1)}{2}}
\alpha\beta
$.
\end{center}

Roughly speaking, $\mathbb{G}_n$ encodes the scalar product of $\mathbb{E}_n$, its orientation\footnote{Encoded by its pseudo-unit $\mathbb{I}_n$.}, and the Grassmann exterior product on $\displaystyle \bigoplus_{k=0}^n \left[\bigwedge^{k} \mathbb{E}_n\right]$, within its associative and distributive algebraic product; moreover, such encoding reveals many connections between algebra and geometry, making new insights possible. 

\section{Notations II}

As we have already said, the Geometric Algebra $\mathbb{G}_n$ is a Euclidean space; we indicate the norm of $X\in\mathbb{G}_n$ with 
symbol~{\color{dgreen}$\bm{|X|}$}$=\sqrt{X\cdot X}$.\index{symboles}{$\vert \cdots \vert$}

In order to emphasize the geometric interpretation of elements in a Geometric algebra, we will use the following nomenclature.

Given an orthonormal ordered basis in the $n$-dimensional Euclidean space $\mathbb{E}_n$ (or an ordered alphabet) $\{\ell_1, \dots , \ell_n\}$, then

\begin{itemize}
	\item elements of $\mathbb{G}_{n \choose 0}\equiv \mathbb{R}$ are called {\bf \color{dgreen} scalars}\index{termes}{Scalars},
	\item elements of $\mathbb{G}_{n \choose 1}\equiv \mathbb{E}_n$ are called {\bf \color{dgreen} vectors}\index{termes}{Vectors},
	\item elements of $\mathbb{G}_{n \choose 2}$ are called {\bf \color{dgreen} bivectors}\index{termes}{Bivector},
	\item elements of $\mathbb{G}_{n \choose k}$ are called {\bf \color{dgreen} $k$-vectors}\index{termes}{$K$-vectors},
	\item elements of $\mathbb{G}_{n \choose n}$ are called {\bf \color{dgreen} pseudo-scalars}\index{termes}{Pseudo-scalars}, and are real multiples of the pseudo-unit $\mathbb{I}_n=\ell_1\cdots\ell_n$.
\end{itemize}

A $k$-vector of the form $\alpha (v_1\wedge \cdots \wedge v_k)$, where $\alpha\in\mathbb{R}$ and each $v_i\in \mathbb{G}_{n \choose 1}$, is called {\bf \color{dgreen} $k$-blade}.\index{termes}{$K$-blade}
Note that in $\mathbb{G}_3$ every bivector is a $2$-blade, while in $\mathbb{G}_4$ the bivector $(\ell_1\ell_2)\ + \ (\ell_3\ell_4)$ is not a $2$-blade. 

In order to limit the use of parentheses, we establish the following precedence rules for operations in $\mathbb{G}_n$, listed below with decreasing rank of precedence\footnote{Sometimes we will also use spacing to stress precedence.}:

\begin{enumerate}
	\item outer product $\wedge$,
	\item product in $\mathbb{G}_n$ among elements of the same grade,
	\item scalar product,
	\item product between a scalar and a $k$-vector (with $k>0$),
	\item sum in $\mathbb{G}_n$.
\end{enumerate}

Thus, for example, $\alpha a \ + \ \beta b$ means $(\alpha a)+(\beta b)$; $\alpha \beta\ x\wedge y$ means $(\alpha \beta)(x\wedge y)$.
However, note also that $\alpha \beta\ x\wedge y = (\alpha x)\wedge (\beta y) =
(\beta x)\wedge (\alpha y)= x\wedge (\alpha\beta y)= \dots$

\section{Inverse of vectors and pseudo-unit in $\mathbb{G}_n$}

In $\mathbb{G}_n$ a vector $v$ is invertible if and only if $v\ne 0$; in this case we have 

\[
v^{-1}
=
\frac{1}{|v|^2} v\ ,
\]

\noindent
where 
{\bf \color{dgreen} $|v|$}$=\sqrt{v\cdot v}=\sqrt{v^2}$. 
In fact $\displaystyle v v^{-1}=\frac{1}{|v|^2} vv=\frac{1}{|v|^2} (v\cdot v + v\wedge v)=1$.

\noindent
In $\mathbb{G}_n$ the pseudo-unit is always invertible, and 

\[
(\mathbb{I}_n)^{-1}
=
(\ell_1\cdots\ell_n)^{-1}
=
\ell_n\cdots\ell_1
=
(-1)^{\frac{n(n-1)}{2}}
\mathbb{I}_n \ .
\]

\section{Pseudo-scalars in $\mathbb{G}_2$ and determinants of $2\times 2$ real matrices}\label{sec:determinants}

In $\mathbb{G}_2$ the notions of bivector, $2$-blade and pseudo-scalar coincide.

Let $\{\ell_1,\ \ell_2\}$ be an ordered orthonormal basis in $\mathbb{E}_2\equiv \mathbb{G}_{2\choose 1}$.

For each $x,y\in \mathbb{G}_{2\choose 1}$, we can write $x=\chi_1 \ell_1\ +\ \chi_2\ell_2$ and $y=\zeta_1 \ell_1\ +\ \zeta_2\ell_2$ (where $\chi_i=x\cdot \ell_i$ and $\zeta_i=y\cdot \ell_i$), then

\[
x\wedge y
=
\left(
	\chi_1 \zeta_2 - \chi_2\zeta_1
\right)
(\ell_1\wedge\ell_2)
=
\det
\left(
\begin{array}{cc}
	\chi_1 & \chi_2 \\
	\zeta_1 & \zeta_2
\end{array}
\right)
\ell_1 \ell_2
=
\det
\left(
\begin{array}{cc}
	\chi_1 & \chi_2 \\
	\zeta_1 & \zeta_2
\end{array}
\right)
\mathbb{I}_2\ ,
\]

and then

\[
(x\wedge y)
(\mathbb{I}_2)^{-1}
=
\det
\left(
\begin{array}{cc}
	\chi_1 & \chi_2 \\
	\zeta_1 & \zeta_2
\end{array}
\right)
=
(x\wedge y)\cdot \mathbb{I}_2\ 
.
\]

Notice that the one-dimensional space $\mathbb{G}_{2\choose 2}$ has the following positive definite symmetric bilinear form:

\[
(x\wedge y)\cdot (w\wedge z)
=
(x\cdot w)(y\cdot z)- (x\cdot z)(y\cdot w)\ ,
\]

where $x,y,w,z\in \mathbb{G}_{2\choose 1}$. 

It will be useful to note the following property, too.
 
\begin{proposition}\label{prop:bivector norm ineq}
If $x,y\in \mathbb{E}_2$, then $|x\wedge y| \le |x| \ |y|$.
\end{proposition}

\textit{\textbf{Proof:}} let $\{\ell_1,\ \ell_2\}$ be an ordered orthonormal basis in $\mathbb{E}_2$; 

let us write $x=\chi_1 \ell_1\ +\ \chi_2\ell_2$ and $y=\zeta_1 \ell_1\ +\ \zeta_2\ell_2$, so that 

\begin{center}
$|x \wedge y|^2 =\left(\chi_1 \zeta_2 - \chi_2\zeta_1 \right)^2$.
\end{center}

The thesis is then achieved verifying the following equivalence,

\begin{center}
$
\left(\chi_1 \zeta_2 - \chi_2\zeta_1 \right)^2 \le (\chi_1^2 +\chi_2^2)(\zeta_1^2+\zeta_2^2)
\Longleftrightarrow
0 \le (\chi_1\zeta_1 + \chi_2\zeta_2)^2 \ . \ \square
$
\end{center}

In this work the symbol {\bf \color{dgreen} $\bm{\square}$}\index{symboles}{$\square$} indicates the end of a proof.

\section{Exterior product of orthogonal vectors in $\mathbb{E}_2$}

\begin{proposition} If $x,y\in \mathbb{E}_2$ are orthogonal, then

\[
x\wedge y
=
\left\{
\begin{array}{rl}
	|x|\ |y| \ \mathbb{I}_2 & \textrm{ if } (x\wedge y)\cdot \mathbb{I}_2 > 0 \\
-	|x|\ |y| \ \mathbb{I}_2 & \textrm{ if } (x\wedge y)\cdot \mathbb{I}_2 < 0
\end{array}
\right.
\]
\end{proposition}

\textit{\textbf{Proof:}} let $\{\ell_1,\ \ell_2\}$ be an ordered orthonormal basis in $\mathbb{E}_2$; $x\wedge y = [(x\wedge y)\cdot\mathbb{I}_2] \mathbb{I}_2$; let us write $x=\chi_1 \ell_1\ +\ \chi_2\ell_2$ and $y=\zeta_1 \ell_1\ +\ \zeta_2\ell_2$, so that $x \wedge y = \left(	\chi_1 \zeta_2 - \chi_2\zeta_1 \right)\mathbb{I}_2$, then

\begin{center}
$
\displaystyle
|\chi_1 \zeta_2 - \chi_2\zeta_1|^2
=
(\chi_1)^2 (\zeta_2)^2+ (\chi_2)^2(\zeta_1)^2 - \chi_1 \zeta_1 \chi_2\zeta_2 - \chi_1 \zeta_1 \chi_2\zeta_2\ ,
$
\end{center}

\noindent
and orthogonality corresponds to the relation $\chi_1\zeta_1=-\chi_2\zeta_2$, so

\begin{center}
$
\displaystyle
|\chi_1 \zeta_2 - \chi_2\zeta_1|^2
=
(\chi_1)^2 (\zeta_2)^2+ (\chi_2)^2(\zeta_1)^2 + (\chi_1)^2 (\zeta_1)^2 + (\chi_2)^2(\zeta_2)^2
=
|x|^2\ |y|^2 \ \square
$
\end{center}

Two ordered bases\footnote{Or, what is the same think, two 2-blades $b_1\wedge b_2$, $c_1\wedge c_2$ in $\mathbb{E}_2$.} $\{b_1,b_2\}$ and $\{c_1,c_2\}$ in $\mathbb{E}_2$ are said to be {\bf \color{dgreen} equi-oriented}\index{termes}{Equi-oriented basis in $\mathbb{E}_2$}\index{termes}{Equi-oriented 2-blades in $\mathbb{E}_2$} if 

\begin{center}
$(b_1\wedge b_2)\cdot (c_1 \wedge c_2)>0$.
\end{center}

\section{Isometric duality between $\mathbb{G}_{3\choose 1}$ and $\mathbb{G}_{3\choose 2}$: the cross product}\label{sec:cross product}

In $\mathbb{G}_3$ the subspaces $\mathbb{G}_{3\choose 1}$ and $\mathbb{G}_{3\choose 2}$ have the same dimension, and are both Euclidean spaces. Moreover, the correspondence 

\[
\mathbb{G}_{3 \choose 1} \ni x  \longmapsto  {\color{dgreen}\bm{x^*}}=x\mathbb{I}_3 \in \mathbb{G}_{3 \choose 2}\ ,
\]\index{symboles}{$x^*$}

\noindent
is an isometry in $\mathbb{G}_3$, whose inverse is the correspondence

\[
\mathbb{G}_{3 \choose 2} 
\ni X \longmapsto {\color{dgreen}\bm{X^{\#}}} = -X \mathbb{I}_3 
\in \mathbb{G}_{3 \choose 1} \ .
\]\index{symboles}{$X^{\#}$}

\noindent
This allows to establish the classical $1:1$ correspondence between a two-dimensional vector subspace of $\mathbb{E}_3$ generated by the two independent vectors 
$a,b\in\mathbb{G}_{3\choose 1}$ (i.e. the bivector $a\wedge b \in \mathbb{G}_{3 \choose 2}$) and the vector $(a\wedge b)^\#\in \mathbb{G}_{3 \choose 1}$ that is orthogonal to that subspace; as a matter of fact

\begin{center}
$ \displaystyle
\big[(a\wedge b)^\# \big]\cdot a 
=
- \big[(a\wedge b)\mathbb{I}_3 \big]\cdot a 
=
-\frac{1}{4}
\big[
(ab-ba)\mathbb{I}_3a + a(ab-ba)\mathbb{I}_3
\big]
=
0\ ,
$
\end{center}

\noindent
as $w\mathbb{I}_3=\mathbb{I}_3w$ for all $w\in \mathbb{G}_{3 \choose 1}$ \Big($\big[(a\wedge b)^\# \big]\cdot b=0$, analogously\Big).
In this sense, the application

\begin{eqnarray*}
\mathbb{G}_{3 \choose 1} \times \mathbb{G}_{3 \choose 1} 
& \longrightarrow &
\mathbb{G}_{3 \choose 1} \\
(a,b) & \longmapsto & (a \wedge b)^\# = - (x\wedge y)\mathbb{I}_3\ ,  
\end{eqnarray*}

corresponds to the classical cross product $a \times b$.

%% file: 03-Geometry.tex
\chapter{Geometric Algebra and Geometry}

\section{Point, lines and planes in $\mathbb{E}_n$}\label{sec:point lines}

In an $n$-dimensional real affine Euclidean space~{\bf \color{dgreen} $\mathbb{A}_n$}\index{symboles}{$\mathbb{A}_n$}, if one fixes a point as the origin, the points in $\mathbb{A}_n$ can be identified with vectors in a Euclidean space $\mathbb{E}_n$. With respect to the foregoing identification, we will talk about  points, lines and planes in a Euclidean space; that is, vectors in $\mathbb{E}_n$  are considered as points of an affine Euclidean space where some point (the corresponding origin $\mathbb{O}$) is fixed somewhere\footnote{See, for instance, ...}.

An {\bf \color{dgreen} oriented (linear) direction}\index{termes}{Oriented (linear) direction in $\mathbb{E}_n$} in $\mathbb{E}_n$ is a nonzero vector.

An {\bf \color{dgreen} oriented (planar) direction}\index{termes}{Oriented (planar) direction in $\mathbb{E}_n$} in $\mathbb{E}_n$ is a nonzero $2$-blade in $\mathbb{G}_n$(the Geometric Algebra associated to $\mathbb{E}_n$).

A {\bf \color{dgreen} line}\index{termes}{Line in $\mathbb{E}_n$}  in $\mathbb{E}_n$ (passing trough point $x_0\in \mathbb{E}_n$ and parallel to the oriented direction $v\in \mathbb{E}_n$) is the set

\[
\big\{
x\in \mathbb{E}_n \ : \ \exists \tau\in \mathbb{R} \ \ \ x=x_0 + \tau v
\big\}\ .
\]

Considering $\mathbb{E}_n\equiv \mathbb{G}_{n \choose 1}$, the same set can be described through vectors and $2$-blades as follows:

\[
\big\{
x\in \mathbb{G}_{n \choose 1} \ : \ (x-x_0) \wedge v =0 
\big\}\ .
\]

Similarly, a {\bf \color{dgreen} plane}\index{termes}{Plane in $\mathbb{E}_n$}  in $\mathbb{E}_n$ (passing trough point $x_0\in \mathbb{E}_n$ and parallel to the vector space generated by the two independent vectors $u,v\in \mathbb{E}_n$) is the set

\[
\big\{
x\in \mathbb{E}_n \ : \ \exists \mu,\ \nu \in \mathbb{R} \  \ x=x_0 + \mu u + \nu v
\big\}\ ,
\]

\noindent
which can also be described through vectors and $3$-vectors in $\mathbb{G}_n$ as follows:

\[
\big\{
x\in \mathbb{G}_{n \choose 1} \ : \ (x-x_0) \wedge u \wedge v =0 
\big\}\ ,
\]

\noindent
characterized by point $x_0 \in \mathbb{E}_n\equiv \mathbb{G}_{n \choose 1}$ and the oriented (plane) direction $u~\wedge~v~\in~\mathbb{G}_{n \choose 2}$.

\section{Oriented intervals and triangles in $\mathbb{E}_n$}

Geometric Algebra is particularly suited to deal with triangles in Euclidean spaces. In particular, the coordinate-free formalism of Geometric Algebra provides analogies between intervals in $\mathbb{R}$ and triangles in $\mathbb{E}_n$. 

Let $\{h_1,\dots , h_n\}$ be an ordered orthonormal basis in the $n$-dimensional Euclidean space $\mathbb{E}_n$. 

An {\bf \color{dgreen} interval}\index{termes}{Interval in $\mathbb{E}_n$}  in $\mathbb{E}_n$ with {\bf \color{dgreen} extremities}\index{termes}{Extremities of an interval} $a,b\in \mathbb{E}_n$ is the set

\[
\big\{
x\in \mathbb{E}_n \ : \  \ x=\alpha a + \beta b \ , \ \alpha + \beta =1 \ , \ \alpha,\beta \ge 0
\big\}\ .
\]

When we want to attribute an orientation to this set, depending on the order of its extremities, we indicate it with the ordered brackets {\bf \color{dgreen} $\bm{[a,b]}$}\index{symboles}{$[a,b]$} and call it  {\bf \color{dgreen} oriented interval}\index{termes}{Oriented interval in $\mathbb{E}_n$} in~$\mathbb{E}_n$. The length of an oriented interval $[a,b]$ is, of course, $|a-b|$.

Analogously, a {\bf \color{dgreen} triangle}\index{termes}{Triangle in $\mathbb{E}_n$}  in $\mathbb{E}_n$ with {\bf \color{dgreen} vertices}\index{termes}{Vertices of a triangle} $a,b,c \in \mathbb{E}_n$ is the set

\[
\big\{
x\in \mathbb{E}_n \ : \  \ x=\alpha a + \beta b + \gamma c\ , \ \alpha + \beta +\gamma =1 \ , \ \alpha,\beta,\gamma \ge 0
\big\}\ .
\]

When we want to attribute an orientation to this set, depending on the order of its vertices, we indicate it within the brackets {\bf \color{dgreen} $\bm{[a,b,c]}$}\index{symboles}{$[a,b,c]$} and call it  {\bf \color{dgreen} oriented triangle}\index{termes}{Oriented triangle in $\mathbb{E}_n$} in~$\mathbb{E}_n$. 

The following oriented triangles

\begin{center}
$
\hfil [b,c,a] \hfil [c,a,b] \hfil
$
\end{center}

\noindent
correspond to the same triangle and have the same orientation with $[a,b,c]$, while the oriented triangles

\begin{center}
$
\hfil [a,c,b] \hfil [c,b,a] \hfil [b,a,c]\hfil
$
\end{center}

\noindent
correspond to the same triangle of $[a,b,c]$, but have opposite orientation with respect the orientation of $[a,b,c]$.

The oriented intervals $[a,b]$, $[b,c]$ and $[c,a]$ are called the {\bf \color{dgreen} sides}\index{termes}{Sides of an oriented triangle} of the oriented triangle $[a,b,c]$.
A {\bf \color{dgreen} diameter}\index{termes}{Diameter of an oriented triangle} of an oriented triangle is a side of maximal length\footnote{An oriented triangle can have more then one diameter, of course.}.

Given an oriented triangle $[a,b,c]$, we consider the bivector

\begin{center}
$
\displaystyle
{\bf \color{dgreen} \left\langle a;b;c\right\rangle}\index{symboles}{$\left\langle a;b;c\right\rangle$}
=
a\wedge b + b\wedge c + c\wedge a
\in \mathbb{G}_{n \choose 2}\ .
$
\end{center}

Such a bivector is a $2$-blade, indeed. For, given the oriented triangle $[a,b,c]$, if we define

\begin{center}
$
\hfil 
{\bf \color{dgreen} \ell_a}\index{symboles}{$\ell_a$}= c-b
\hfil 
{\bf \color{dgreen} \ell_b}\index{symboles}{$\ell_b$}= a-c
\hfil 
{\bf \color{dgreen} \ell_c}\index{symboles}{$\ell_c$}= b-a\ ,
\hfil 
$
\end{center}

then 

\begin{center}
$
\displaystyle
\left\langle a;b;c\right\rangle
=
\ell_a \wedge \ell_b
=
\ell_b \wedge \ell_c
=
\ell_c \wedge \ell_a \ .
$
\end{center}

Notice that 

\begin{center}
$
\displaystyle
\left\langle a;b;c\right\rangle
=
\left\langle b;c;a\right\rangle
=
\left\langle c;a;b\right\rangle
=
-\left\langle a;c;b\right\rangle
=
-\left\langle c;b;a\right\rangle
=
-\left\langle b;a;c\right\rangle \ ,
$
\end{center}

that is, $\left\langle a;b;c\right\rangle$ change sign if we change the orientation of $[a,b,c]$. 

The {\bf \color{dgreen} area}\index{termes}{Area of a  triangle} of a triangle whose vertices are $a,b,c$ is $\displaystyle \frac{1}{2}\Big|\left\langle a;b;c\right\rangle\Big|$. Moreover, such area can also be expressed using only the scalar product

\begin{center}
$
\displaystyle 
\frac{1}{2}
\Big|\left\langle a;b;c\right\rangle\Big|
=
\frac{1}{2}
\sqrt{(\ell_a\wedge \ell_b)\cdot(\ell_a\wedge \ell_b)}
=
\frac{1}{2}
\sqrt{|\ell_a|^2 |\ell_b|^2- (\ell_a\cdot \ell_b)^2}\ .
$
\end{center}

A triangle, whose area is zero, is called {\bf \color{dgreen} degenerate}\index{termes}{Degenerate triangle}.
\bigskip

Let us now express the bivector $\left\langle a;b;c\right\rangle \in \mathbb{G}_{n \choose 2}$ coordinatewise.

If  $\{h_1,\dots , h_n\}$ is an ordered orthonormal basis in the $n$-dimensional Euclidean space $\mathbb{E}_n$, and 

\begin{center}
$
\displaystyle
\hfil 
a= \sum_{j=1}^n \alpha_j h_j
\hfil 
b= \sum_{j=1}^n \beta_j h_j
\hfil 
c= \sum_{j=1}^n \gamma_j h_j\ ,
\hfil 
$
\end{center}

then 

\[
\begin{array}{l}
\phantom{=}
	\left\langle a;b;c\right\rangle
=
a\wedge b + b\wedge c + c\wedge a = \\ \\
\displaystyle
=
\sum_{1\le j< k \le n}
\left[
 \det
  \left(
   \begin{array}{cc}
		 \alpha_j & \alpha_{k} \\ 
		 \beta_j & \beta_{k} 
	 \end{array}
	\right)
	+
\det
  \left(
   \begin{array}{cc}
		 \beta_j & \beta_{k} \\ 
		 \gamma_j & \gamma_{k} 
	 \end{array}
	\right)
	+
\det
  \left(
   \begin{array}{cc}
		 \gamma_j & \gamma_{k} \\
		 \alpha_j & \alpha_{k} 
	 \end{array}
	\right)
	+
\right]
h_j \wedge h_{k} = \\ \\
\displaystyle
=
(b-a)\wedge (c-a) 
=
\sum_{1\le j< k \le n}
 \det
  \left(
   \begin{array}{cc}
		 \beta_j -\alpha_j & \beta_{k}-\alpha_{k} \\
		 \gamma_j - \beta_j & \gamma_{k}- \beta_{k} 
	 \end{array}
	\right)
h_j \wedge h_{k} \ .
\end{array}
\]

In particular,

\begin{center}
$
\displaystyle 
\frac{1}{2}
\Big|\left\langle a;b;c\right\rangle\Big|
=
\frac{1}{2}
\sqrt{
\sum_{1\le j< k \le n}
\left[
 \det
  \left(
   \begin{array}{cc}
		 \beta_j -\alpha_j & \beta_{k}-\alpha_{k} \\
		 \gamma_j - \beta_j & \gamma_{k}- \beta_{k} 
	 \end{array}
	\right)
\right]^2
}\ ,
$
\end{center}

since $\displaystyle \big\{h_j\wedge h_{k}\big\}_{1\le j< k \le n}$ is an orthonormal basis in $\mathbb{G}_{n \choose 2}$.

\section{Reflections in $\mathbb{E}_n$ and mirror vertices in plane triangles}\label{sec:mirror vertices}

Invertible vectors (that is, linear directions) are useful to represent mirror points with respect to those directions.

Let us consider a point $x \in \mathbb{E}_n$ and a direction $v \in \mathbb{E}_n$; if $x$ and $v$ are linearly independent ($x\wedge v \ne 0$), then the mirror image of $x$ with respect to the line passing through $0$ and $v$ is the point 

\psset{unit=0.8cm}

\begin{center}
\begin{figure}[!h]
\begin{pspicture}(-4.5,1.9)(7,6.1)
\psdots(1.96,1.82)(5.68,4.1)(3.95,5.71) 
\uput[d](1.96,1.82){$0$}
\uput[r](5.68,4.1){$x$}
\uput[u](3.95,5.71){$vxv^{-1}$}
\uput[r](3.2,3.16){$v$}
\psline{->}(1.96,1.82)(3.2,3.16)
\psline{->}(1.96,1.82)(5.68,4.1)
\psline[linestyle=dotted]{->}(1.96,1.82)(3.95,5.71)
\psline[linestyle=dotted](1.24,1.04)(6.2,6.4)
\psline[linestyle=dotted](5.68,4.1)(3.95,5.71)
\end{pspicture}
\end{figure}
\end{center}

\[
\begin{array}{rl}
	vxv^{-1}
	&
	\displaystyle
	=
	(v\cdot x \ + \ v \wedge x)v^{-1}
	=
	(v\cdot x \ - \ x \wedge v)v^{-1}
	= \\ \\
	&
	\displaystyle
	=
	\big[v\cdot x \ - \ (xv- \ x\cdot v)\big] v^{-1}
	=
	(2 v\cdot x \ - \ xv)v^{-1}
	=
	2\frac{x\cdot v}{|v|^2}v - x
\end{array}
\] 

\noindent
(see~\cite{Lou2001} at pag.13 for further details). Note that the foregoing formula works even when $x\wedge v=0$ (that is, $x=\chi v$ for some $\chi \in \mathbb{R}$).

Let us consider a nondegenerate oriented triangle $[a,b,c]$ in $\mathbb{E}_2$ (that is, a plane triangle), then each of its (oriented) sides $[a,b]$, $[b,c]$ and $[c,a]$ determines a direction ($\ell_c$, $\ell_a$ and $\ell_b$ respectively). So we can consider the mirror image {\bf \color{dgreen} $\bm{x'}$}\index{symboles}{$x'$} of each vertex $x$ of $[a,b,c]$ with respect to the line passing through its two adjacent vertices. We have that  

\[
\begin{array}{l}
\displaystyle
a'
=
c + \ell_a	\ell_b \ell_a^{-1} = c + 2 \frac{\ell_a\cdot \ell_b}{|\ell_a|^2}\ell_a -\ell_b 
=
-\left[
a 
+ 2 b \frac{\ell_b\cdot \ell_a}{|\ell_a|^2} 
+ 2 c \frac{\ell_c\cdot \ell_a}{|\ell_a|^2}
\right]\ , \\ \\
\displaystyle
b'
=
a + \ell_b	\ell_c \ell_b^{-1} = a + 2 \frac{\ell_b\cdot \ell_c}{|\ell_b|^2}\ell_b -\ell_c 
=
-\left[
b 
+ 2 c \frac{\ell_c\cdot \ell_b}{|\ell_b|^2} 
+ 2 a \frac{\ell_a\cdot \ell_b}{|\ell_b|^2}
\right]\ , \\ \\
\displaystyle
c'
=
b + \ell_c	\ell_a \ell_c^{-1} = b + 2 \frac{\ell_c\cdot \ell_a}{|\ell_c|^2}\ell_c -\ell_a 
=
-\left[
c 
+ 2 a \frac{\ell_a\cdot \ell_c}{|\ell_c|^2} 
+ 2 b \frac{\ell_b\cdot \ell_c}{|\ell_c|^2}
\right]\ ,
\end{array}
\]

and we call them {\bf \color{dgreen} mirror vertices}\index{termes}{Mirror vertices of a triangle} of the oriented\footnote{However, they do not depend on the triangle's orientation.} triangle $[a,b,c]$. 
\bigskip

Given a nondegenerate oriented plane triangle $[a,b,c]$, it will be useful to define the following point and two directions

$
\displaystyle
\hfill
{\bf \color{dgreen} \bar{a}}\index{symboles}{$\bar{a}$}
=
\frac{1}{2}(a'+a)\ ,
\hfill
{\bf \color{dgreen} u_a}\index{symboles}{$u_a$}
=
\frac{1}{2}(a'-a)\ ,
\hfill
{\bf \color{dgreen} v_a}\index{symboles}{$v_a$}
=
c-\bar{a}\ .
\hfill
$

The above definitions can be generalized to any vertex. Indeed, if $x$ is a vertex of a nondegenerate oriented plane triangle $[x,x_+,x_-]$ one can define 

\[
\begin{array}{rl}
\displaystyle
x'
& \displaystyle 
=
x_- + \ell_{x}	\ell_{x_+} \ell_{x}^{-1} = x_- + 2 \frac{\ell_{x}\cdot \ell_{x_+}}{|\ell_x|^2}\ell_x -\ell_{x_+}
= \\
& \displaystyle 
=
-\left[
x 
+ 2 (x_+) \frac{\ell_{x_+}\cdot \ell_x}{|\ell_x|^2} 
+ 2 (x_-) \frac{\ell_{x_-}\cdot \ell_x}{|\ell_x|^2}
\right]\ , \\ \\
\displaystyle
\bar{x}
& \displaystyle 
=
\frac{1}{2}(x'+x)
\ \ ,\ \
u_x
=
\frac{1}{2}(x'-x)
\ \ , \ \
v_x
=
x_- - \bar{x}\ ,
\end{array}
\]

where $\ell_x= (x_- -x_+)$, $\ell_{x_+}=(x- x_-)$ and $\ell_{x_-}=(x_+ - x)$.
So we have that

\begin{center}
$
\displaystyle
\hfill
x = \bar{x} - u_x
\hfill
x_+ = \bar{x} + (v_x -\ell_x)
\hfill
x_- = \bar{x} + v_x
\hfill
x'= \bar{x}+ u_x\ , 
\hfill
$
\end{center}

and we can state the following proposition.

\begin{proposition}
Let $x$ be a vertex of the nondegenerate oriented plane triangle $[x,x_+,x_-]$, then

\begin{enumerate}
	\item 
	$
  2 u_x \wedge \ell_x
  =	
	2 \left\langle x;x_+;x_-\right\rangle
	=
	\left\langle x;x_+;x_-\right\rangle
	-
	\left\langle x';x_+;x_-\right\rangle\ ;
	$
	\item $u_x \cdot \ell_x=0$ \big(so that $u_x\wedge \ell_x= \pm |u_x| \ |\ell_x|\ \mathbb{I}_2$\big)\ ;
	\item $v_x \wedge \ell_x=0$ \big(so there exists $\tau\in\mathbb{R}$ such that $v_x=\tau\ell_x$ \big)\ .
\end{enumerate}

\end{proposition}

\textbf{\textit{Proof of 2.}} We have that $\displaystyle u_x=\frac{1}{2}\big(x' -x\big) = \frac{1}{2}\big(-\ell_{x_+} + \ell_x \ell_{x_+}\ell_x^{-1}\big)$, and

\begin{center}
$
\begin{array}{l}
	\displaystyle
	\phantom{=}
	u_x \cdot \ell_x 
	=
	\frac{1}{2}
	\big(u_x\ell_x + \ell_x u_x\big)
	=
	\frac{1}{4}
	\big(
	-\ell_{x_+}\ell_x +\ell_x \ell_{x_+} \ell_x^{-1} \ell_x
	-\ell_x \ell_{x_+}+ \ell_x\ell_x \ell_{x_+}\ell_x^{-1}
	\big)
	=0\ .
\end{array}
$
\end{center}

\textbf{\textit{Proof of 3.}} We have that $\displaystyle v_x= x_- - \frac{1}{2}\big(x'+x \big)= -\frac{1}{2}\big(\ell_{x_+}+ \ell_x \ell_{x_+}\ell_x^{-1}\big)$, and

\begin{center}
$
\begin{array}{l}
	\displaystyle
	\phantom{=}
	\kern-10pt
	v_x \wedge \ell_x 
	=
	\frac{1}{2}
	\big(v_x\ell_x - \ell_x v_x\big)
	=
	\frac{1}{4}
	\big(
	-\ell_{x_+}\ell_x -\ell_x \ell_{x_+} \ell_x^{-1} \ell_x
	+\ell_x \ell_{x_+}+ \ell_x\ell_x \ell_{x_+}\ell_x^{-1}
	\big)
	=0 .\ \square
\end{array}
$
\end{center}

A mirror vertex $x'$ of the oriented triangle $[x,x_+,x_-]$ is said to be {\bf \color{dgreen} balanced}\index{termes}{Balanced mirror vertex of a triangle} if $\bar{x}\in [x_+,x_-]$. In particular, if $[x_+,x_-]$ is a diameter of that triangle, then $x'$ is balanced. This implies that every triangle has at least a balanced mirror vertex.
Owing to the foregoing proposition, balanced mirror vertices can be characterized through lengths.

\begin{proposition}\label{pro:mirror vertex}
A mirror vertex $x'$ of the oriented triangle $[x,x_+,x_-]$ is balanced if and only if $|\ell_x|= |\ell_x-v_x|+ |v_x|$.
\end{proposition}


In the two following figures, the mirror vertices $a'$ and $b'$ are balanced, while the mirror vertex $c'$ is not.

\psset{unit=0.85cm}

\begin{figure}[!h]
\
\hspace{50pt}
\begin{pspicture}(-3,0)(1,6)
\psdots(0,0)(1,3)(-1,3) 
\pspolygon(0,0)(1,3)(-1,3)
\psdots(0,6)
\psdots(-2.6,1.8)
\psdots(0,3)
\psdots(-0.8,2.4)
\uput[d](0,0){$a$}
\uput[r](1,3){$b$}
\uput[l](-1,3){$c$}
\uput[u](0,6){$a'$}
\uput[dl](-2.6,1.8){$b'$}
\uput[ur](0,3){$\bar{a}$}
\uput[dl](-0.7,2.4){$\bar{b}$}
\psline[linestyle=dashed](0,0)(0,6) 
\psline[linestyle=dashed](1,3)(-2.6,1.8) 
\pspolygon[linestyle=dotted](1,3)(0,6)(-1,3)
\pspolygon[linestyle=dotted](-1,3)(-2.6,1.8)(0,0)
\end{pspicture}
\hspace{60pt}
\begin{pspicture}(-3,-3)(3,1)
\psdots(0,0)(3,1)(-3,1) 
\pspolygon(0,0)(3,1)(-3,1)
\psdots(-1.8,-2.6) 
\psdots(-2.4,-0.8) 
\uput[d](0,0){$a$}
\uput[r](3,1){$b$}
\uput[l](-3,1){$c$}
\uput[d](-1.8,-2.6){$c'$}
\uput[l](-2.4,-0.8){$\bar{c}$}
\psline[linestyle=dashed](-3,1)(-1.8,-2.6) 
\psline[linestyle=dotted](0,0)(-1.8,-2.6)(3,1)
\end{pspicture}
\end{figure}

\vfill

%% file: 04-SmoothCurves.tex
\chapter{Smooth curves}\label{cha:smooth curves}

In the following sections we describe some classical approximation algorithms for smooth curves in $\mathbb{E}_n$ because of their analogies with our approximation Algorithms~\ref{eq:AlgI} and~\ref{eq:AlgII} for smooth surfaces in $\mathbb{E}_n$.

\section{Approximation through inscribed mean vectors}

Let $I \subseteq \mathbb{R}$ be an open interval. Let $\mathbb{E}_n$ be a $n$-dimensional Euclidean space, and $\{h_1,\dots ,h_n\}$ an orthonormal basis in $\mathbb{E}_n$.

A continuous function $c:I\to \mathbb{E}_n$ will be simply called a {\bf \color{dgreen} curve}.\index{termes}{Curve in $\mathbb{E}_n$}

In this chapter we define $n$ real functions $\gamma_j:I\to \mathbb{R}$ as the $n$ components $\gamma_j=c\cdot h_j$; so, we have that $\forall \tau \in I$ $\displaystyle c(\tau)=\sum_{j=1}^n \gamma_j(\tau) h_j$.

A vector $v\in \mathbb{E}_n$ is said to be {\bf \color{dgreen} inscribed}\index{termes}{Inscribed vector} in the curve $c$ if there exist $\alpha,\beta\in I$, such that $\alpha\ne \beta$ and $v=c(\beta)-c(\alpha)$; in this case the vector $\displaystyle \frac{1}{\beta-\alpha}[c(\beta)-c(\alpha)]$ is called {\bf \color{dgreen} inscribed mean vector}\index{termes}{Inscribed mean vector} in $c$.

A curve $c:I\to \mathbb{E}_n$ is said to be {\bf \color{dgreen} smooth}\index{termes}{Smooth curve} if

\begin{itemize}
	\item each $\gamma_j$ has continuous second derivative $\ddot{\gamma}_j$ on $I$, and
	\item there exists $\delta >0$ such that  $\displaystyle \sup_{\tau\in I}|\ddot{\gamma}_j(\tau)|\le \delta <\infty$ ($\delta$ does not depend on $j=1,\dots, n$).
\end{itemize}

For a smooth curve the following estimate holds for each $\tau, \tau+\epsilon \in I$

\begin{equation}\label{eq:estimate curve}
\big|
\gamma_j(\tau+\epsilon) -\gamma_j(\tau)-\dot{\gamma}_j(\tau) \epsilon
\big|
\le
\frac{\delta}{2} \epsilon^2\ , 
\end{equation}

\noindent
where $\dot{\gamma_j}$ is, of course, the first derivative of $\gamma_j$. Following the Landau notation, we can rewrite the relation~(\ref{eq:estimate curve}) as follows

\[
\gamma_j(\tau+\epsilon) -\gamma_j(\tau)= \dot{\gamma}_j(\tau) \epsilon + O(\epsilon^2)\ .
\]

If $c$ is a smooth curve, we denote $\displaystyle {\bf \color{dgreen} \dot{c}(\tau)}\index{symboles}{$\dot{c}(\tau)$}=\sum_{j=1}^n \dot{\gamma_j}(\tau)\ h_j$.

\begin{proposition}\label{prop: approxim dot c}
If $c:I\to \mathbb{E}_n$ is a smooth curve and $\chi \in I$, then the vector $\dot{c}(\chi)$ is the limit of the inscribed mean vectors $\displaystyle \frac{1}{\beta-\alpha}\big[c(\beta)-c(\alpha)\big]$ as $(\alpha,\beta) \to (\chi,\chi)$ in $\mathbb{R}^2$, that is

\[
\lim_{
\begin{array}{c}
\scriptstyle 	(\alpha,\beta) \to (\chi,\chi) \\
\scriptstyle	\alpha \ne \beta
\end{array}
}
\frac{1}{\beta-\alpha}\big[c(\beta)-c(\alpha)\big]
=
\dot{c}(\chi) \ .
\]

\end{proposition}

The proof of Proposition~\ref{prop: approxim dot c} is routine. Nonetheless we prove it, just because the proof we provide here is a $1$-dimensional version of the proof of Theorem~\ref{thm:II}.
\bigskip

\textit{\textbf{Proof of Proposition~\ref{prop: approxim dot c}:}} let $\alpha \ne \beta$, then

\[
\begin{array}{l}
	\displaystyle
	\frac{\gamma_j(\beta)-\gamma_j(\alpha)}{\beta - \alpha}
	=
	\frac{\gamma_j(\beta)-\gamma_j\left(\frac{\alpha+\beta}{2}\right) -\left[\gamma_j(\alpha)-\gamma_j\left(\frac{\alpha+\beta}{2}\right)\right]}{\beta - \alpha}
	= \\ \\
	\displaystyle
	=
	\frac{\dot{\gamma_j}\left(\frac{\alpha+\beta}{2}\right)\left(\beta - \frac{\alpha+\beta}{2}\right) + O\left(\left(\beta - \frac{\alpha+\beta}{2}\right)^2\right)}{\beta-\alpha} -\frac{\dot{\gamma_j}\left(\frac{\alpha+\beta}{2}\right)\left(\alpha - \frac{\alpha + \beta}{2}\right)+O\left(\left(\alpha - \frac{\alpha + \beta}{2}\right)^2\right)}{\beta - \alpha}= \\ \\
	\displaystyle
	=
\dot{\gamma_j}\left(\frac{\alpha+\beta}{2}\right) + O(\beta -\alpha)\ .
\end{array}
\]

So, 

\[
\begin{array}{l}
\displaystyle
	\left| \frac{1}{\beta -\alpha}\big[c(\beta)-c(\alpha)\big] - \dot{c}(\chi)\right|
	=
	\left| \sum_{j=1}^n \left[\dot{\gamma_j}\left(\frac{\alpha+\beta}{2}\right) 
	- \dot{\gamma_j}(\chi) 
	+ O(\beta -\alpha)\right] h_j\right| = \\ \\
\displaystyle
\le 
O(\beta -\alpha) + 	
\sum_{j=1}^n \left| \dot{\gamma_j}\left(\frac{\alpha+\beta}{2}\right) 
	- \dot{\gamma_j}(\chi) \right|
\end{array}
\]

that goes to $0$ as $(\alpha,\beta)\to(\chi,\chi)$ because each $\gamma_j$ is $C^2(I)$\ $\square$

\section{Geometric interpretation of the direction $\dot{c}(\chi)$}

If $c$ is a smooth curve and $\dot{c}(\chi) \ne 0$, then $c$ is locally injective and thus, if $\alpha$ and $\beta$ are sufficiently close to $\chi \in I$, every inscribed mean vector $\displaystyle \frac{1}{\beta -\alpha}\big[c(\beta)-c(\alpha)\big]$ is a direction, and Proposition~\ref{prop: approxim dot c} has the following geometric interpretation:

direction $\dot{c}(\chi)$ is the limit of the inscribed mean directions $\displaystyle \frac{1}{\beta -\alpha}\big[c(\beta)-c(\alpha)\big]$ as $(\alpha,\beta)\to (\chi,\chi)$.

It is well known that if $c$ is smooth, but $\dot{c}(\chi)= 0$, then the foregoing interpretation may be false, even if $c$ is locally injective. A classical example is the cusp in~$\mathbb{E}_2$, $c(\tau)= \tau^2 h_1 + \tau^3 h_2$, and $\chi=0$.

\section{Estimates of the length of a smooth curve}

If the curve $c:I\to \mathbb{E}_n$ is smooth and $[\alpha,\beta]\subset I$, then the following algorithm 

\begin{equation}\label{eq: Algorithm 0}
\sum_{i=0}^k\big|c(\beta_i)-c(\alpha_i)\big|
\end{equation}

can estimate the integral\footnote{That we call {\bf \color{dgreen} length}\index{termes}{Length of a smooth curve} of the curve $c:[\alpha,\beta]\to\mathbb{E}_n$, when $c$ is injective on $[\alpha,\beta]$.} $\displaystyle \int_\alpha^\beta \big|\dot{c}(\tau)\big|\ d\tau$, where 

\begin{itemize}
	\item $\alpha_0=\alpha$,
	\item $\alpha_i < \alpha_{i+1}=\beta_i$ (for $i=0,\dots, k-1$), and
	\item $\beta_{k-1}< \beta_k =\beta$;
\end{itemize}

that is, $\big\{[\alpha_i,\beta_i]\big\}_{i=0}^k=\Pi$ is a partition of $[\alpha.\beta]$ with contiguous nonoverlapping intervals; in this sense Algorithm~(\ref{eq: Algorithm 0}) can also be written

\[
\sum_{[\alpha_i,\beta_i]\in \Pi}\big|c(\beta_i)-c(\alpha_i)\big|\ .
\]

More precisely, Algorithm~(\ref{eq: Algorithm 0}) converges to $\displaystyle \int_\alpha^\beta \big|\dot{c}(\tau)\big|\ d\tau$ when the maximal length of intervals in the partition $\Pi$, $\displaystyle \max_{[\alpha_i,\beta_i]\in\Pi}|\beta_i-\alpha_i|$, goes to zero, as we can see from the following elementary estimates:

\begin{eqnarray*}
& & 
	\left|
		\sum_{[\alpha_i,\beta_i]\in \Pi}\big|c(\beta_i)-c(\alpha_i)\big| 
		-
		\int_\alpha^\beta \big|\dot{c}(\tau)\big|\ d\tau
	\right|
	= \\ 
&	= &
	\left|
		\sum_{[\alpha_i,\beta_i]\in \Pi}
		\left[
		\big|c(\beta_i)-c(\alpha_i)\big| 
		-
		\int_{\alpha_i}^{\beta_i} \big|\dot{c}(\tau)\big|\ d\tau
		\right]
	\right|
	= \\
& \le & 
	\sum_{[\alpha_i,\beta_i]\in \Pi}
		\left|
		\big|c(\beta_i)-c(\alpha_i)\big| 
		-
		\int_{\alpha_i}^{\beta_i} \big|\dot{c}(\tau)\big|\ d\tau
		\right|
	= \\
\end{eqnarray*}

\begin{eqnarray*}
&	= &
	\sum_{[\alpha_i,\beta_i]\in \Pi}
		\left|
		\int_{\alpha_i}^{\beta_i}
		\left[
		\frac{\big|c(\beta_i)-c(\alpha_i)\big|}{|\beta_i-\alpha_i|}
		-
		 \big|\dot{c}(\tau)\big|
		\right]\ d\tau
		\right|
	= \\
&	\le &
	\sum_{[\alpha_i,\beta_i]\in \Pi}
		\int_{\alpha_i}^{\beta_i}
		\left|
		\frac{\big|c(\beta_i)-c(\alpha_i)\big|}{|\beta_i-\alpha_i|}
		-
		 \big|\dot{c}(\tau)\big|
		\right|\ d\tau = (*)\ .
\end{eqnarray*}

As $\Big||v|-|w|\Big|\le |v-w|$ for each $v,w\in \mathbb{E}_n$, we have that 

\begin{eqnarray*}
	(*) & \le &
	\sum_{[\alpha_i,\beta_i]\in \Pi}
		\int_{\alpha_i}^{\beta_i}
		\left|
		\frac{1}{\beta_i-\alpha_i}
		\big[c(\beta_i)-c(\alpha_i)\big]
		-
		\dot{c}(\tau)
		\right|\ d\tau
	= \\
&	= &
	\sum_{[\alpha_i,\beta_i]\in \Pi}
		\int_{\alpha_i}^{\beta_i}
		\left|
		\frac{1}{\beta_i-\alpha_i}
		\big[c(\beta_i)-c(\alpha_i)-\dot{c}(\alpha_i)(\beta_i-\alpha_i)\big]
		+\big[
		\dot{c}(\alpha_i)
		 -
		\dot{c}(\tau)
		\big]
		\right|\ d\tau
	= \\
&	\le &
	\sum_{[\alpha_i,\beta_i]\in \Pi}
		\int_{\alpha_i}^{\beta_i}
		\left[
		\left|
		\frac{1}{\beta_i-\alpha_i}
		\big[c(\beta_i)-c(\alpha_i)-\dot{c}(\alpha_i)(\beta_i-\alpha_i)\big]
		\right|
		+
		\left|
		\dot{c}(\alpha_i)
		 -
		\dot{c}(\tau)
		\right|
		\right]\ d\tau \ .
\end{eqnarray*}

Then, owing to estimate~(\ref{eq:estimate curve}), we have that 

\[
\begin{array}{l}
\displaystyle
	\Big|
	c(\beta_i)-c(\alpha_i)-\dot{c}(\alpha_i)(\beta_i-\alpha_i)
	\Big|
	= \\ \\
\displaystyle 	
=
	\left|
	\sum_{j=1}^n
	\big[\gamma_j(\beta_i)-\gamma_j(\alpha_i)-\dot{\gamma_j}(\alpha_i)(\beta_i-\alpha_i)\big]h_j
	\right|
	= \\ \\
	\displaystyle 	
	\le
	\sum_{j=1}^n
	\Big|
	\gamma_j(\beta_i)-\gamma_j(\alpha_i)-\dot{\gamma_j}(\alpha_i)(\beta_i-\alpha_i)
\Big|
	\le n \delta (\beta_i-\alpha_i)^2\ .
\end{array}
\]

So we can conclude that

\[
\begin{array}{l}
	\displaystyle 
	\left|
		\sum_{[\alpha_i,\beta_i]\in \Pi}\big|c(\beta_i)-c(\alpha_i)\big| 
		-
		\int_\alpha^\beta \big|\dot{c}(\tau)\big|\ d\tau
	\right|
	= \\ \\
	\displaystyle 
	\le
	n\delta
	\sum_{[\alpha_i,\beta_i]\in \Pi}
	 (\beta_i-\alpha_i)^2
	+
	\sum_{[\alpha_i,\beta_i]\in \Pi}
	\int_{\alpha_i}^{\beta_i}
		\left|
		\dot{c}(\alpha_i)
		 -
		\dot{c}(\tau)
		\right|\ d\tau \ .
\end{array}
\]

That goes to zero as $\displaystyle \max_{[\alpha_i,\beta_i]\in \Pi} |\beta_i-\alpha_i|\longrightarrow 0$, since $c$ is smooth.
\bigskip

Here we provided, as we did before for the proof of Proposition~\ref{prop: approxim dot c}, many details of well-known estimates, just because of their analogies with the estimates for the area of a smooth surface\footnote{See section~\ref{sec:area}.}.

%% file: 05-SmoothSurfaces.tex
\chapter{Smooth surfaces}

\section{Surfaces, inscribed balanced mean bivectors}\label{sec:surfaces}

Let $\Omega$ be an open set in $\mathbb{E}_2$.
A continuous function $s:\Omega \to \mathbb{E}_n$ will be simply called a {\bf \color{dgreen} surface}.\index{termes}{Surface in $\mathbb{E}_n$}

If $\{h_1,\dots , h_n\}$ is an ordered basis in the $n$-dimensional Euclidean space $\mathbb{E}_n$ and $s$ is a surface, we define $n$ real functions $\sigma_j:\Omega \to \mathbb{R}$ as the components $\sigma_j = s \cdot h_j$; so, we have that $\displaystyle\forall x \in \Omega\ \ s(x)=\sum_{j=1}^n \sigma_j(x) h_j$, and $s$ is a surface if and only if each component $\sigma_j$ is continuous.
If $\{\ell_1, \ell_2\}$ is an ordered orthonormal basis in $\mathbb{E}_2$, we can indicate each $x\in\mathbb{E}_2$ as
$\chi_1 \ell_1 + \chi_2 \ell_2 $, where $\chi_1=x\cdot \ell_1$ and $\chi_2=x\cdot \ell_2$.

\begin{exam}\label{exa:cylinder}
A  circular right cylinder of radius $\rho$ is a surface. Thus, it corresponds to the following function $s:\mathbb{E}_2\to \mathbb{E}_3$, 

\begin{center}
$
s(x)
=
s(\chi_1 \ell_1 + \chi_2 \ell_2)=\rho \cos(\chi_1)h_1 + \rho \sin(\chi_1)h_2 + \chi_2h_3
$,
\end{center}

\noindent
where $\{\ell_1,\ell_2\}$ is an ordered orthonormal basis in $\mathbb{E}_2$, $\{h_1,h_2,h_3\}$ is an ordered orthonormal basis in $\mathbb{E}_3$.
\end{exam}

A bivector $V\in\mathbb{G}_{n \choose 2}$ is said to be {\bf \color{dgreen} inscribed}\index{termes}{Inscribed bivector in a surface} in a surface $s$ if there exists a nondegenerate ordered plane triangle $[a,b,c]$ contained in $\Omega$, such that

\begin{center}
$
\displaystyle
V
=
\big\langle 
s(a);s(b);s(c)
\big\rangle
=
s(a)\wedge s(b)
\ +\ 
s(b)\wedge s(c)
\ +\ 
s(c)\wedge s(a)\ .
$
\end{center}

In this case, the ordered triangle $\big[s(a),s(b),s(c)\big]$ is said to be {\bf \color{dgreen} inscribed}\index{termes}{Inscribed ordered triangle on a surface} on the surface $s$, and the following bivector

\begin{center}
$
\displaystyle
\frac{1}{
\left\langle 
a;b;c
\right\rangle
\cdot \mathbb{I}_2
}
\big\langle 
s(a);s(b);s(c)
\big\rangle
$
\end{center}

\noindent
is called {\bf \color{dgreen} inscribed mean bivector}\index{termes}{Inscribed mean bivector in a surface} in $s$.

\noindent
A plane triangle $\Delta$ of vertices $\{a,b,c\}$ is said to be {\bf \color{dgreen} balanced}\index{termes}{Balanced triangle in a open set of $\mathbb{E}_2$} in the open set~$\Omega$ if 

\begin{itemize}
	\item $\Delta \subset \Omega$,
	\item there exists a balanced mirror vertex $x'\in \{a',b',c'\}$ of $\Delta$ such that the plane triangle of vertices\footnote{See Section~\ref{sec:mirror vertices} for the notations.} $\{x',x_+,x_-\}$ is contained in $\Omega$.
\end{itemize}

An inscribed triangle $\big[s(a),s(b),s(c)\big]$ is said to be {\bf \color{dgreen} balanced}\index{termes}{Inscribed balanced ordered triangle on a surface} on $s$ if the plane triangle $[a,b,c]$ is balanced in $\Omega$. Note that, since $\Omega$ is an open set, every sufficiently small\footnote{Small with respect to its diameter.} inscribed triangle $\big[s(a),s(b),s(c)\big]$ is balanced on $s$.

If $\big[s(a),s(b),s(c)\big]$ is an inscribed  balanced ordered triangle on $s$, where $a'$ is a balanced mirror vertex of $[a,b,c]$ (with respect to vertex $a$) such that $[a',b,c]$ is in $\Omega$, then the bivector

\begin{center}
$
\displaystyle
\big[s(a')-s(a)\big]
\wedge
\big[s(c)-s(b)\big]
$
\end{center}

\noindent
is called {\bf \color{dgreen} inscribed balanced bivector}\index{termes}{Inscribed balanced bivector in a surface} in $s$;
in this case the bivector

\begin{center}
$
\displaystyle
\frac{1}{2\left\langle a;b;c\right\rangle\cdot \mathbb{I}_2}
\big[s(a')-s(a)\big]
\wedge
\big[s(c)-s(b)\big]
$
\end{center}

\noindent
is called {\bf \color{dgreen} inscribed balanced mean bivector}\index{termes}{Inscribed balanced mean bivector in a surface} in $s$.

An inscribed balanced bivector can be written in different ways

\[
\begin{array}{l}
\displaystyle
\big[s(a')-s(a)\big]
\wedge
\big[s(c)-s(b)\big]
= \\ \\
=
s(a)\wedge s(b) 
\ + \	
s(b)\wedge s(a') 
\ + \	
s(a')\wedge s(c) 
\ + \	
s(c)\wedge s(a)
= \\ \\
\displaystyle
=
\big\langle 
s(a);s(b);s(c)
\big\rangle
+
\big\langle 
s(c);s(b);s(a')
\big\rangle
= 
\big\langle 
s(a);s(b);s(c)
\big\rangle
-
\big\langle 
s(a');s(b);s(c)
\big\rangle\ .
\end{array}
\]

It will also be useful to write an inscribed balanced bivector coordinatewise with respect to an ordered orthonormal basis $\{h_1, \dots , h_n\}$ in~$\mathbb{E}_n$; that is, if $\displaystyle\forall x \in \Omega$ $\displaystyle s(x)=\sum_{j=1}^n \sigma_j(x) h_j$, then

\[
\begin{array}{l}
\displaystyle
\phantom{=}
\big[s(a')-s(a)\big]
\wedge
\big[s(c)-s(b)\big]
= \\ \\
\displaystyle
=
\left\{\sum_{j=1}^n \big[\sigma_j(a')-\sigma_j(a)\big] h_j\right\}	
\wedge
\left\{\sum_{k=1}^n \big[\sigma_k(c)-\sigma_k(b)\big] h_k\right\}	
= \\ \\
\displaystyle
=
\sum_{1\le j<k\le n}
\Big\{
\big[\sigma_j{\scriptstyle(a')} \kern-2pt - \kern-2pt \sigma_j{\scriptstyle(a)}\big]
\big[\sigma_{k}{\scriptstyle(c)} \kern-2pt - \kern-2pt \sigma_{k}{\scriptstyle(b)}\big] 
\kern-2pt - \kern-2pt
\big[\sigma_j{\scriptstyle(c)} \kern-2pt - \kern-2pt \sigma_j{\scriptstyle(b)}\big]
\big[\sigma_{k}{\scriptstyle(a')} \kern-2pt - \kern-2pt \sigma_{k}{\scriptstyle(a)}\big]
\Big\}
h_j\wedge h_{k}	\ .
\end{array}
\]

Let us now define $n \choose 2$ transformations ${\color{dgreen} \bm{s_{j,k}}}:\Omega \to \mathbb{E}_2$,\index{symboles}{$s_{j,k}$}

\begin{center}
$
\displaystyle
s_{j,k}(x)
=
\sigma_j(x) \ell_1 \ + \ \sigma_{k}(x) \ell_2\ ;
$
\end{center}

then we can rewrite each component

\[
\begin{array}{l}
\displaystyle
	\big[\sigma_j(a')-\sigma_j(a)\big]\big[\sigma_{k}(c)-\sigma_{k}(b)\big] 
	-
	\big[\sigma_j(c)-\sigma_j(b)\big]\big[\sigma_{k}(a')-\sigma_{k}(a)\big]
	= \\ \\
	\displaystyle
	=
	\Big\{
	\big[s_{j,k}(a')-s_{j,k}(a)\big]
	\wedge
	\big[s_{j,k}(c)-s_{j,k}(b)\big]
	\Big\}
	\cdot
	\mathbb{I}_2\ ,
\end{array}
\]

\noindent
so that an inscribed balanced bivector can be written as

\[
\big[s(a')-s(a)\big]
\wedge
\big[s(c)-s(b)\big]
=
\kern-10pt
\sum_{1\le j<k\le n} 
\kern-5pt
\Big\{
	\big\{
	[s_{j,k}(a')-s_{j,k}(a)]
	\wedge
	[s_{j,k}(c)-s_{j,k}(b)]
	\big\}
	\cdot
	\mathbb{I}_2
\Big\}
h_j\wedge h_{k}	 \ .
\] 

An inscribed balanced mean bivector can be written in other ways, too; in particular \big(as $\left\langle a;b;c\right\rangle = \left\langle c;b;a'\right\rangle$\big) it corresponds to the following mean of inscribed mean bivectors:

\begin{center}
$
\begin{array}{l}
	\displaystyle
	\frac{1}{2\left\langle a;b;c\right\rangle\cdot \mathbb{I}_2}
	\big[s(a')-s(a)\big]
	\wedge
	\big[s(c)-s(b)\big]
	= \\ \\
	\displaystyle
	=
	\frac{1}{2}
	\left\{
	\frac{1}{\left\langle a;b;c\right\rangle\cdot \mathbb{I}_2}
	\big\langle s(a);s(b);s(c)\big\rangle
	+
	\frac{1}{\left\langle c;b;a'\right\rangle\cdot \mathbb{I}_2}
	\big\langle s(c);s(b);s(a')\big\rangle
\right\} \ .
\end{array}
$
\end{center}

In the case of a surface in space ($s:\Omega \to \mathbb{E}_3$) the foregoing mean can be seen as the mean of the vectors orthogonal to the planes secant the surfaces at points $s(a), s(b), s(c)$ and $s(c), s(b), s(a')$ respectively. However, the most interesting representation of an inscribed balanced bivector on $s$ is the following:

\begin{center}
$
\begin{array}{l}
	\displaystyle
	\frac{1}{2\left\langle a;b;c\right\rangle\cdot \mathbb{I}_2}
	\big[s(a')-s(a)\big]
	\wedge
	\big[s(c)-s(b)\big]
	= \\ \\
	\displaystyle
	=
	\frac{1}{\big[ \left\langle a;b;c\right\rangle - \left\langle a';b;c\right\rangle\big]\cdot \mathbb{I}_2}
	\Big[
	\big\langle s(a);s(b);s(c)\big\rangle
	-
	\big\langle s(a');s(b);s(c)\big\rangle
 \Big] \ ;
\end{array}
$
\end{center}

\noindent
indeed, the previous expression plays for surfaces (in Theorem~\ref{thm:II}) the same role as the classical expression (for the inscribed mean vector)

\begin{center}
$
\displaystyle
\frac{1}{\beta -\alpha}
\big[
c(\beta) - c(\alpha)
\big]
$
\end{center}

\noindent
plays for curves (in Proposition~\ref{prop: approxim dot c}).

\section{Notations III}

Let $\Omega\subseteq \mathbb{E}_2$ be open. Let us give some differential notations for a sufficiently regular function $\psi:\Omega \to \mathbb{R}$. Let $w \in \mathbb{E}_2$ 

\[
\begin{array}{rcll}
{\bf \color{dgreen} \partial_w\psi(x)}\index{symboles}{$\partial_w\psi(x)$} 
& = &
\displaystyle \lim_{\epsilon\to 0} \frac{1}{\epsilon} \big[\psi(x+\epsilon w)-\psi(x)\big] \in \mathbb{R} & \\ \\
{\bf \color{dgreen} \nabla \psi(x)}\index{symboles}{$\nabla \psi(x)$} 
& = &
\partial_{\ell_1}\psi(x) \ell_1 
+
\partial_{\ell_2}\psi(x) \ell_2
\in \mathbb{E}_2 
& 
\textrm{ is the gradient vector,}
\\ \\
{\bf \color{dgreen} \bm{H}_\psi(x)}\index{symboles}{$\bm{H}_\psi(x)$} 
& = &
\left(
\begin{array}{cc}
	\partial_{\ell_1}\partial_{\ell_1}\psi(x)
	&
	\partial_{\ell_2}\partial_{\ell_1}\psi(x) \\ \\
	\partial_{\ell_1}\partial_{\ell_2}\psi(x)
	&
	\partial_{\ell_2}\partial_{\ell_2}\psi(x)
\end{array}
\right) \in \mathbb{R}^{2 \times 2}
& 
\textrm{ is the Hessian matrix,}
\end{array}
\]

\noindent
whose real eigenvalues are indicated as {\bf \color{dgreen} $\lambda_{i,\psi(x)}$}\index{symboles}{$\lambda_{i,\psi(x)}$} 
(with $i=1,2$), corresponding to the ordered orthonormal basis of eigenvectors $\{\ell_{1,\psi(x)}, \ell_{2,\psi(x)}\}$ equioriented with a fixed ordered orthonormal basis $\{\ell_1, \ell_2\}$ in $\mathbb{E}_2\equiv \mathbb{G}_{2 \choose 1}$, that is $\ell_{1,\psi(x)}\wedge \ell_{2,\psi(x)} = \ell_1 \wedge \ell_2= \mathbb{I}_2$.

\section{Smooth functions, transformations and surfaces}

Let $\Omega\subseteq \mathbb{E}_2$ be open, a (real) {\bf \color{dgreen} function} $\psi:\Omega \to \mathbb{R}$ is said to be {\bf \color{dgreen} smooth}\index{termes}{Smooth (real) function} if

\begin{itemize}
	\item $\psi$ has second-order derivatives, that are continuous on $\Omega$ \big(is $C^2(\Omega)$\big),
	\item $\displaystyle \sup_{x\in \Omega} \max\big\{ |\lambda_{1,\psi(x)}|, |\lambda_{2,\psi(x)}| \big\}< +\infty$.
\end{itemize}

For a smooth function $\psi$ the following estimate\footnote{That estimate is analogue to estimate~(\ref{eq:estimate curve}), and can be obtained using the Taylor formula with integral remainder.} holds for each $x,y\in \Omega$ such that the interval $[x,y]$ is contained in $\Omega$

\begin{equation}\label{eq:estimate function}
\big|
\psi(y)-\psi(x) - \nabla\psi(x) \cdot (y-x)
\big|
\le
\lambda
|y-x|^2\ ,
\end{equation}

where $\displaystyle \lambda =  \sup_{z\in [x,y]} \max\big\{ |\lambda_{1,\psi(z)}|, |\lambda_{2,\psi(z)}| \big\}$. 

Following the Landau notation, we can rewrite relation~(\ref{eq:estimate function}) as follows\footnote{Here, for the sake of simplicity, we can suppose that the set $\Omega$ is also convex.}:

\begin{equation}\label{eq:estimate function Landau}
\forall x,y\in \Omega
\hspace{30pt}
\psi(y)-\psi(x) - \nabla\psi(x) \cdot (y-x)
=
O(|y-x|^2)\ .
\end{equation}

In particular, if $z,w\in \mathbb{E}_2$ are such that the interval $[z-w,w+w]$ is contained in $\Omega$, we can also write

\[
\psi(z+w)-\psi(z-w) = 2 \nabla\psi(z) \cdot w \ +\ O(|w|^2)\ .
\]

Let $\Omega\subseteq \mathbb{E}_2$ be open, let $\{\ell_1,\ell_2\}$ be an ordered orthonormal basis in $\mathbb{E}_2$, and $\{h_1,\dots h_n\}$ be an ordered orthonormal basis in $\mathbb{E}_n$,

\begin{itemize}
	\item a {\bf \color{dgreen} transformation} $f:\Omega \to \mathbb{E}_2$ is called {\bf \color{dgreen} smooth}\index{termes}{Smooth plane transformation} if each of its components $\phi_i=f\cdot e_i$ is a smooth (real) function (i=1,2);
	\item a {\bf \color{dgreen} surface} $s:\Omega \to \mathbb{E}_n$ is called {\bf \color{dgreen} smooth}\index{termes}{Smooth surface} if each of its components $\sigma_j=s\cdot h_i$ is a smooth (real) function (j=1,\dots , n).
\end{itemize}

\begin{exam}
The cylinder of example~\ref{exa:cylinder} is a smooth surface. As a matter of fact, if $\{\ell_1,\ell_2\}$ is an ordered orthonormal basis in $\mathbb{E}_2$ and $x=\chi_1\ell_1 + \chi_2\ell_2$, then

\begin{center}
$
\hfil
\sigma_1(x)= \rho\cos\chi_1\ , 
\hfil
\sigma_2(x)= \rho\sin\chi_1\ , 
\hfil
\sigma_3(x)= \chi_2\ ,
\hfil
$
\end{center}

and $\displaystyle \sup_{x\in \mathbb{E}_2} \max_{1\le j \le 3} \max_{1\le i \le 2} |\lambda_{i,\sigma_j(x)}|=\rho$.
\end{exam}

%% file: 06-MainResults.tex
\chapter{Main results}

\section{The Jacobian determinant}

\begin{theorem}\label{thm:I}
Let $\{\ell_1,\ell_2\}$ be an ordered orthonormal  basis in the Euclidean space $\mathbb{E}_2$;
let $\Omega \subseteq \mathbb{E}_2$ be open;  let $f:\Omega \to \mathbb{E}_2$ be a smooth plane transformation, then $\forall x \in \Omega$ we have that

\begin{equation}\label{eq:Alg0}
	\lim_{
	\begin{array}{c}
		\scriptstyle (a,b,c)\to (x,x,x) \\
		\scriptstyle 	a\wedge b +b\wedge c + c\wedge a \ne 0
	\end{array}
	}
	\kern-20pt
	\frac{\Big\{\big[f(d_{(a,b,c)})-f(a)\big]\wedge \big[f(c)-f(b)\big]\Big\}\cdot \mathbb{I}_2}{2\big(a\wedge b +b\wedge c + c\wedge a\big)\cdot \mathbb{I}_2}
	=
	\big(\nabla\phi_1(x)\wedge \nabla\phi_2(x)\big)\cdot \mathbb{I}_2
\end{equation}

\noindent
where 

\begin{enumerate}
	\item $\displaystyle d_{(a,b,c)}
	=
	a'
	=
-\left[
a 
+ 2 b \frac{\ell_b\cdot \ell_a}{|\ell_a|^2} 
+ 2 c \frac{\ell_c\cdot \ell_a}{|\ell_a|^2}
\right]$ is the mirror vertex of vertex $a$, in the oriented plane triangle~$[a,b,c]$,
	\item $a'$ is balanced,
\end{enumerate}

\noindent
and

\begin{itemize}
	\item $\mathbb{I}_2= \ell_1 \ell_2=\ell_1\wedge \ell_2$ is the pseudo-unit in the Geometric Algebra $\mathbb{G}_2$ associated to the oriented Euclidean space $\mathbb{E}_2$,
	\item $\wedge$ is the outer product in $\mathbb{G}_2$, and $\cdot$ is the scalar product in $\mathbb{G}_2$,
	\item $\phi_i = f\cdot \ell_i$ (with $i=1,2$),
	\item the limit $(a,b,c)\to (x,x,x)$ is taken in the product topology of $\mathbb{E}_2 \times \mathbb{E}_2 \times \mathbb{E}_2$.
\end{itemize}

\end{theorem}

\begin{rem}
Coordinatewise, the transformation $f$ can be seen as the real transformation $\bm{\Phi}:\tilde{\Omega}\to \mathbb{R}^2$, where

\begin{itemize}
	\item $\mathbb{R}^2 \supseteq \tilde{\Omega} \ni (\chi_1,\chi_2) \Longleftrightarrow \chi_1 \ell_1 \ +\ \chi_2 \ell_2 \in \Omega \subseteq \mathbb{E}_2$;
	\item $
		\displaystyle 
			\bm{\Phi}(\chi_1,\chi_2)
			 = 
	\big(f(\chi_1 \ell_1 \ +\ \chi_2 \ell_2)\cdot \ell_1\ , \ f(\chi_1 \ell_1 \ +\ \chi_2 \ell_2)\cdot \ell_2\big)$
		
		$\displaystyle 
		\phantom{	\bm{\Phi}(\chi_1,\chi_2)}
		=
		\big(
		\phi_1(\chi_1 \ell_1 \ +\ \chi_2 \ell_2)\ , \ \phi_2(\chi_1 \ell_1 \ +\ \chi_2 \ell_2)
		\big)
		$;
\end{itemize}

then, we have that

\begin{center}
$
\displaystyle
\left(
\begin{array}{cc}
\partial_{\ell_1}\phi_1(x) & \partial_{\ell_2}\phi_1(x) \\ \\ 	
\partial_{\ell_1}\phi_2(x) & \partial_{\ell_2}\phi_2(x) 
\end{array}
\right)
=
\frac{\partial (\phi_1,\phi_2)}{\partial(\chi_1,\chi_2)}
$ 
\end{center}

is the Jacobian matrix of transformation $\bm{\Phi}$, whose determinant\footnote{See Section~\ref{sec:determinants}.} is

\[
\big(\nabla\phi_1(x)\wedge \nabla\phi_2(x)\big)\cdot \mathbb{I}_2
=
\det
\frac{\partial (\phi_1,\phi_2)}{\partial(\chi_1,\chi_2)}
\ .
\]
\end{rem}

\begin{rem}
It is possible to obtain an analogue result for transformations within $k$-di\-men\-sional Euclidean spaces. However, such result will be treated in other works, were we will apply it to construct $k$-dimensional Stieltjes-like measures in $\mathbb{E}_n$.
\end{rem}

\begin{rem}
The approximating ratio in the foregoing theorem can be rewritten as a kind of incremental ratio for the function $f:\Omega\subseteq\mathbb{E}_2\to\mathbb{E}_2$

\[
\frac{\Big\{\big[f(a')-f(a)\big]\wedge \big[f(c)-f(b)\big]\Big\}\cdot \mathbb{I}_2}{2\big(a\wedge b +b\wedge c + c\wedge a\big)\cdot \mathbb{I}_2}
=
\frac
{\Big[
	\big\langle f(a);f(b);f(c)\big\rangle
	-
	\big\langle f(a');f(b);f(c)\big\rangle
 \Big]\kern-2pt \cdot \kern-2pt \mathbb{I}_2
}
{\big[ \left\langle a;b;c\right\rangle - \left\langle a';b;c\right\rangle\big]\kern-2pt \cdot \kern-2pt \mathbb{I}_2}
\]
where $\left\langle x;y;z\right\rangle=x\wedge y + y\wedge z + z\wedge x$. This strengthens the analogy between the derivative of a function of one real variable and the Jacobian determinant of a transformation of two real variables.
\end{rem}

\begin{rem}
The coordinatewise writing of the previous approximating ratio is more elaborate and lengthy. Let us define $\tilde{\phi}_i(\chi_1,\chi_2)=\phi_i(\chi_1\ell_1+\chi_2\ell_2)$, then $\tilde{\phi}_i:\tilde{\Omega}\to\mathbb{R}$. If we denote $a=\alpha_1\ell_1+\alpha_2\ell_2$, $a'=\alpha'_1\ell_1+\alpha'_2\ell_2$, $b=\beta_1\ell_1+\beta_2\ell_2$, $c=\gamma_1\ell_1+\gamma_2\ell_2$, the thesis of Theorem~\ref{thm:I} becomes

\[
\frac
{\det
\left(
\begin{array}{cc}
	\tilde{\phi}_1(\alpha'_1,\alpha'_2)-\tilde{\phi}_1(\alpha_1,\alpha_2) 
	&
	\tilde{\phi}_2(\alpha'_1,\alpha'_2)-\tilde{\phi}_2(\alpha_1,\alpha_2) \\
	\tilde{\phi}_1(\gamma_1,\gamma_2)-\tilde{\phi}_1(\beta_1,\beta_2) 
	&
	\tilde{\phi}_2(\gamma_1,\gamma_2)-\tilde{\phi}_2(\beta_1,\beta_2)
\end{array}
\right)}
{2
\det
\left(
\begin{array}{cc}
	\beta_1-\alpha_1 & \beta_2-\alpha_2 \\
	\gamma_1-\beta_1 & \gamma_2 -\beta_2
\end{array}
\right)}
\longrightarrow
\det
\frac{\partial (\tilde{\phi}_1,\tilde{\phi}_2)}{\partial(\chi_1,\chi_2)}
\]

as $(\alpha_1,\alpha_2,\beta_1,\beta_2,\gamma_1,\gamma_2)\longrightarrow(\chi_1,\chi_2,\chi_1,\chi_2,\chi_1,\chi_2)$ in $\mathbb{R}^6$, where

\[
\alpha'_i
=
-
\left[
\alpha_i
+
2\beta_i
\frac{\scriptstyle (\alpha_1 -\gamma_1)(\gamma_1-\beta_1)+(\alpha_2 -\gamma_2)(\gamma_2-\beta_2)}{\scriptstyle (\gamma_1-\beta_1)^2+(\gamma_2-\beta_2)^2}
+2\gamma_i
\frac{\scriptstyle (\beta_1 -\alpha_1)(\gamma_1-\beta_1)+(\beta_2 -\alpha_2)(\gamma_2-\beta_2)}{\scriptstyle (\gamma_1-\beta_1)^2+(\gamma_2-\beta_2)^2}
\right]\ .
\]

The comparison between the above Cartesian expressions and the Geometric Algebraic ones should suggest why we prefer the Clifford coordinate-free language.
\end{rem}

\textit{\textbf{Proof of Theorem~\ref{thm:I}.}}

As $\mathbb{E}_2$ is locally convex, every sufficiently small triangle is balanced in $\Omega$. Let us write the inscribed balanced bivector with respect to the ordered basis $\{\ell_1, \ell_2\}$,

\begin{center}
$
\begin{array}{l}
	\displaystyle 
	\phantom{=}
	\big[f(a')-f(a)\big]\wedge \big[f(c)-f(b)\big]= \\ 
	\displaystyle
	=
	\Big\{
		\big[\phi_1{\scriptstyle(a')}-\phi_1{\scriptstyle(a)}\big]\ell_1 
			+ \kern-3pt
		\big[\phi_2{\scriptstyle(a')}-\phi_2{\scriptstyle(a)}\big]\ell_2
	\Big\}
	\wedge
	\Big\{
		\big[\phi_1{\scriptstyle(c)}-\phi_1{\scriptstyle(b)}\big]\ell_1 
			+ \kern-3pt
		\big[\phi_2{\scriptstyle(c)}-\phi_2{\scriptstyle(b)}\big]\ell_2
	\Big\}
	= \\ 
	\displaystyle
	=
	\Big\{
	\big[\phi_1{\scriptstyle(a')}-\phi_1{\scriptstyle(a)}\big]
	\big[\phi_2{\scriptstyle(c)}-\phi_2{\scriptstyle(b)}\big]
	-
	\big[\phi_2{\scriptstyle(a')}-\phi_2{\scriptstyle(a)}\big]
	\big[\phi_1{\scriptstyle(c)}-\phi_1{\scriptstyle(b)}\big]
	\Big\}
	\mathbb{I}_2
\end{array}
$
\end{center}

putting $\bar{a}=\frac{1}{2}(a'+a)$, $u_a=\frac{1}{2}(a'-a)$, $v_a=c-\bar{a}$ and $\ell_a=c-b$ we have that 

$a'=\bar{a}+u_a$, $a=\bar{a}-u_a$, $c=\bar{a}+v_a$ and $b=\bar{a}-(\ell_a-v_a)$.
As $\mathbb{E}_2$ is locally convex, there exists in $\Omega$ a convex open neighborhood of $x$ where we can use estimate~(\ref{eq:estimate function Landau}) 

\begin{eqnarray*}
\phi_i(a')-\phi_i(a)
& = &
\phi_i(\bar{a}+u_a)-\phi_i(\bar{a}-u_a)
=
\nabla\phi_i(\bar{a})\cdot (2u_a)\ + \ O\big(|u_a|^2\big) \\ \\
	\phi_i(c)-\phi_i(b)
&	= &
	\phi_i(\bar{a}+v_a)-\phi_i\big(\bar{a}-(\ell_a-v_a)\big)
	= \\
&	= &
	\phi_i(\bar{a}+v_a)-\phi_i(\bar{a})-\Big[\phi_i\big(\bar{a}-(\ell_a-v_a)\big)-\phi_i(\bar{a})\Big]
	= \\
&	= &
	\nabla\phi_i(\bar{a})\cdot v_a \ + \ O\big(|v_a|^2\big) +\nabla\phi_i(\bar{a})\cdot (\ell_a-v_a)\ +\ O\big(|\ell_a-v_a|^2\big)
	= \\
&	= &
	\nabla\phi_i(\bar{a})\cdot \ell_a \ + \ O\big(|v_a|^2\big) \ +\ O\big(|\ell_a-v_a|^2\big) \ .
\end{eqnarray*}

Then,

\[
\begin{array}{l}
	\displaystyle 
	\phantom{=}
	\Big\{\big[f(a')-f(a)\big]\wedge \big[f(c)-f(b)\big]\Big\}\cdot \mathbb{I}_2= 
	\\ \\
	\displaystyle
	=
	\big[\phi_1{\scriptstyle(a')}-\phi_1{\scriptstyle(a)}\big]
	\big[\phi_2{\scriptstyle(c)}-\phi_2{\scriptstyle(b)}\big]
	-
	\big[\phi_2{\scriptstyle(a')}-\phi_2{\scriptstyle(a)}\big]
	\big[\phi_1{\scriptstyle(c)}-\phi_1{\scriptstyle(b)}\big]
	=	\\ \\
	\displaystyle
	=
		\Big[\nabla\phi_1(\bar{a})\cdot (2u_a)\ + \ O\big(|u_a|^2\big)\Big]
		\Big[\nabla\phi_2(\bar{a})\cdot \ell_a \ + \ O\big(|v_a|^2\big) \ +\ O\big(|\ell_a-v_a|^2\big)\Big] + \\
	\displaystyle
	\phantom{=}
		-
		\Big[\nabla\phi_2(\bar{a})\cdot (2u_a)\ + \ O\big(|u_a|^2\big)\Big]
		\Big[	\nabla\phi_1(\bar{a})\cdot \ell_a \ + \ O\big(|v_a|^2\big) \ +\ O\big(|\ell_a-v_a|^2\big)\Big]= \\ \\
	\displaystyle
	=
	\Big[\nabla\phi_1(\bar{a})\cdot (2u_a)\Big]
	\Big[\nabla\phi_2(\bar{a})\cdot \ell_a\Big]
	-
	\Big[\nabla\phi_2(\bar{a})\cdot (2u_a)\Big]
	\Big[	\nabla\phi_1(\bar{a})\cdot \ell_a\Big] +\\
	\phantom{=}
	+ 
	\Big[\nabla\phi_1(\bar{a})\cdot (2u_a)\Big]
	\Big[O\big(|v_a|^2\big) \ +\ O\big(|\ell_a-v_a|^2\big)\Big]
	+O\big(|u_a|^2\big)\Big[\nabla\phi_2(\bar{a})\cdot \ell_a\Big] + \\
	\phantom{=}
	+ 
	\Big[\nabla\phi_2(\bar{a})\cdot (2u_a)\Big]
	\Big[O\big(|v_a|^2\big) \ +\ O\big(|\ell_a-v_a|^2\big)\Big]
	+O\big(|u_a|^2\big)\Big[\nabla\phi_1(\bar{a})\cdot \ell_a\Big] + \\
	\phantom{=}
	+
	O\big(|u_a|^2\big)\Big[O\big(|v_a|^2\big) \ +\ O\big(|\ell_a-v_a|^2\big)\Big] \ .
\end{array}	
\]

The first difference is a scalar product between bivectors in $\mathbb{G}_2$ (i.e. pseudo-scalars)

\begin{eqnarray*}
	& &
		\Big[\nabla\phi_1(\bar{a})\cdot (2u_a)\Big]
		\Big[\nabla\phi_2(\bar{a})\cdot \ell_a\Big]
		-
		\Big[\nabla\phi_2(\bar{a})\cdot (2u_a)\Big]
		\Big[	\nabla\phi_1(\bar{a})\cdot \ell_a\Big]
	= \\
	& = &
	\Big(\nabla\phi_1(\bar{a})\wedge \nabla\phi_2(\bar{a})\Big)
	\cdot
	\big((2u_a)\wedge \ell_a\big)\ .
\end{eqnarray*}

Since $a'$ is a mirror vertex, then $u_a$ is a direction orthogonal to the direction $\ell_a$, and we have that 

\[
	\phantom{=}
	\Big(\nabla\phi_1(\bar{a})\wedge \nabla\phi_2(\bar{a})\Big)
	\cdot
	\big((2u_a)\wedge \ell_a\big)\\
	= \pm 2 |u_a| |\ell_a| 
	\Big(\nabla\phi_1(\bar{a})\wedge \nabla\phi_2(\bar{a})\Big)
	\cdot
	\mathbb{I}_2\ ,
\]

\noindent
the sign depending on whether or not the ordered basis $\{ u_a, \ell_a \}$ is equi-oriented with the ordered orthononormal basis $\{\ell_1, \ell_2\}$.

Now, we simply observe that 

\[
2\big(a\wedge b +b\wedge c + c\wedge a\big)
=
(2u_a)\wedge \ell_a
=
\pm 2 |u_a| |\ell_a|\  \mathbb{I}_2\ ,
\]

and so we can write

\begin{equation}\label{eq:estimate transf}
\begin{array}{l}
	\displaystyle
	\phantom{=}
		\left|
		\frac{\Big\{\big[f(a')-f(a)\big]\wedge \big[f(c)-f(b)\big]\Big\}\cdot \mathbb{I}_2}{2\big(a\wedge b +b\wedge c + c\wedge a\big)\cdot \mathbb{I}_2}
		-
		\Big(\nabla\phi_1(x)\wedge \nabla\phi_2(x)\Big)
		\cdot
		\mathbb{I}_2
		\right|
	= \\ \\
	\displaystyle
	\le 
		\left|
		\Big(\nabla\phi_1(\bar{a})\wedge \nabla\phi_2(\bar{a})\Big)
		\cdot
		\mathbb{I}_2
		-
		\Big(\nabla\phi_1(x)\wedge \nabla\phi_2(x)\Big)
		\cdot
		\mathbb{I}_2
		\right| + \\
	\displaystyle
	\phantom{\le}
	+
	\big|\nabla\phi_1(\bar{a})\big|
	\left|
		\frac{O\big(|v_a|^2\big)}{|\ell_a|} + \frac{O\big(|\ell_a-v_a|^2\big)}{|\ell_a|}
	\right|
		+
	\big|\nabla\phi_2(\bar{a})\big|\ \Big|O\big(|u_a|\big)\Big| + \\
	\displaystyle
	\phantom{\le\ }
	+
	\big|\nabla\phi_2(\bar{a})\big|
	\left|
		\frac{O\big(|v_a|^2\big)}{|\ell_a|} + \frac{O\big(|\ell_a-v_a|^2\big)}{|\ell_a|}
	\right|
		+ 
	\big|\nabla\phi_1(\bar{a})\big| \ \Big|O\big(|u_a|\big)\Big| + \\ 
	\displaystyle
	\phantom{\le \ \ \ }
	+
	O\big(|u_a|\big)\
	\left|
		\frac{O\big(|v_a|^2\big)}{|\ell_a|} + \frac{O\big(|\ell_a-v_a|^2\big)}{|\ell_a|} 
	\right|
\end{array}
\end{equation}

By Cauchy-Schwarz inequality, triangular inequality and Proposition~\ref{prop:bivector norm ineq}, we have that $\forall t,w,y,z\in \mathbb{E}_2$

\begin{center}
$
\begin{array}{l}
	\phantom{=}
	\big|(t\wedge w)\cdot \mathbb{I}_2 -(y\wedge z)\cdot \mathbb{I}_2\big|
	=
	|t\wedge w -y\wedge z|
	=
	|t\wedge w - t\wedge z + t\wedge z - y\wedge z| = \\
	\le
	\big|t\wedge (w-z)\big| + \big|(t-y) \wedge z\big|
	\le
	|t| |w-z| + |t-y||z|
\end{array}
$
\end{center}

so

\begin{center}
$
\begin{array}{l}
	\displaystyle
	\phantom{\le}
	\left|
		\Big(\nabla\phi_1(\bar{a})\wedge \nabla\phi_2(\bar{a})\Big)
		\cdot
		\mathbb{I}_2
		-
		\Big(\nabla\phi_1(x)\wedge \nabla\phi_2(x)\Big)
		\cdot
		\mathbb{I}_2
	\right| = \\
	\displaystyle
	\le
	\Big|\nabla\phi_1(\bar{a})\Big| \Big|\nabla\phi_2(\bar{a})-\nabla\phi_2(x)\Big|
	+
	\Big|\nabla\phi_1(\bar{a})-\nabla\phi_1(x)\Big| \Big|\nabla\phi_2(x)\Big|
\end{array}
$
\end{center}

Moreover, the mirror vertex $a'$ is balanced, and then (by Proposition~\ref{pro:mirror vertex}) we have that 

\begin{equation}\label{eq:balanced estimate}
\max\left\{\frac{|v_a|}{|\ell_a|}, \frac{|\ell_a-v_a|}{|\ell_a|}\right\}\le 1\ .
\end{equation}

As $f$ is smooth, the theorem is proved; in fact, 

\[
\bar{a}\longrightarrow x
\ \ \ \textrm{ and } \ \ \ 
|u_a|,\ |v_a|,\ |\ell_a-v_a|, |\ell_a| \longrightarrow 0\ ,
\]

as $(a,b,c)\longrightarrow (x,x,x)\ .\ \square$
\bigskip

\begin{rem}\label{rem:relaxing hypothesis}
As we have warned in the introduction\footnote{See Section~\ref{sec:warnings}} some hypotheses in the foregoing theorem can be weakened. For instance, we can relax hypothesis~($\textit{2.}$) by imposing that 

\begin{center}
$
\displaystyle
\max\left\{\frac{|v_a|}{|\ell_a|}, \frac{|\ell_a-v_a|}{|\ell_a|}\right\}
$
\end{center}

\noindent
be simply bounded, instead of supposing $a'$ being balanced \big(that is equivalent to~(\ref{eq:balanced estimate})\big). We will show in section~\ref{sec:not balanced bivector} that even when $d=d_{(a,b,c)}$ is not a mirror vertex of~$[a,b,c]$, it is possible to choose suitable oriented triangles $[d,c,b]$, adjacent to the converging triangles $[a,b,c]$, such that~(\ref{eq:Alg0}) and~(\ref{eq:AlgI}) hold.
\end{rem}

\section{The tangent bivector}

The following theorem is to smooth surfaces as Proposition~\ref{prop: approxim dot c} is to smooth curves.

\begin{theorem}\label{thm:II}
Let $\{\ell_1,\ell_2\}$ be an ordered orthonormal  basis in the Euclidean space $\mathbb{E}_2$; let $\{h_1,\dots , h_n\}$ be an ordered orthonormal  basis in the $n$-dimensional Euclidean space $\mathbb{E}_n$;
let $\Omega \subseteq \mathbb{E}_2$ be open, and $s:\Omega \to \mathbb{E}_n$ be a smooth surface, 
then $\forall x \in \Omega$ we have that

\begin{equation}\label{eq:AlgI}
\begin{array}{c}
\displaystyle
\ \kern-10pt
	\lim_{
	\begin{array}{c}
		\scriptstyle (a,b,c)\to (x,x,x) \\
		\scriptstyle 	a\wedge b +b\wedge c + c\wedge a \ne 0
	\end{array}
	}
	\kern-12pt
	\frac{1}{\big[ \big\langle a;b;c\big\rangle - \big\langle d_{(a,b,c)};b;c\big\rangle\big] \kern-3pt \cdot \kern-1pt \mathbb{I}_2}
		\Big[
		\kern-3pt
		\big\langle s(a);s(b);s(c)\big\rangle
		-
		\big\langle s(d_{(a,b,c)});s(b);s(c)\big\rangle
		\kern-3pt
		\Big]
	\kern-3pt = \\
	\displaystyle
	=
	\partial_{\ell_1} s(x) \wedge \partial_{\ell_2} s(x) \ ,
\end{array}
\end{equation}

where

\begin{itemize}
	\item $\displaystyle d_{(a,b,c)}
	=a'=
-\left[
a 
+ 2 b \frac{\ell_b\cdot \ell_a}{|\ell_a|^2} 
+ 2 c \frac{\ell_c\cdot \ell_a}{|\ell_a|^2}
\right]$ is a balanced mirror vertex of the oriented plane triangle~$[a,b,c]$,
	\item $\mathbb{I}_2= \ell_1 \ell_2=\ell_1\wedge \ell_2$ is the pseudo-unit in the Geometric Algebra $\mathbb{G}_2$ associated to the oriented Euclidean space $\mathbb{E}_2$,
	\item $\wedge$ is the outer product in $\mathbb{G}_k$, and $\cdot$ is the scalar product in $\mathbb{G}_k$ (with $k=2$ or $k=n$),
	\item $
\displaystyle
	{\bf \color{dgreen} \partial_{\ell_i}s(x)}\index{symboles}{$\partial_{\ell_i}s(x)$}
=
\sum_{j=1}^n \partial_{\ell_i} \sigma_j(x)\ h_j
$ (with $i=1,2$), where $\sigma_j = s \cdot h_j$  (with $j=1,\dots , n$),
	\item the limit $(a,b,c)\to (x,x,x)$ is taken in the product topology of $\mathbb{E}_2 \times \mathbb{E}_2 \times \mathbb{E}_2$.
\end{itemize}

\end{theorem}

\textit{\textbf{Proof of Theorem~\ref{thm:II}.}}

The proof is just a coordinatewise application of Theorem~\ref{thm:I}. Let us rewrite

\begin{eqnarray}
	& & 
	\phantom{=}
	\displaystyle
	\big\langle s(a);s(b);s(c)\big\rangle
	-
	\big\langle s(a');s(b);s(c)\big\rangle
	=
	\big[s(a')-s(a)\big] 
	\wedge
	\big[s(c)-s(b)\big] 
	= \nonumber \\
	& & 
	\displaystyle
	=
	\left\{
	\sum_{j=1}^n
	\big[
	\sigma_j(a')-\sigma_j(a)
	\big]
	h_j
	\right\}
	\wedge
	\left\{
	\sum_{k=1}^n
	\big[
	\sigma_k(c)-\sigma_k(b)
	\big]
h_k
\right\} 
= \nonumber \\
	& & 
	\displaystyle
	=
	\kern-7pt
	\sum_{1\le j< k \le n}
	\kern-7pt
	\Big\{
	\kern-3pt
	\big[
	\sigma_j(a')-\sigma_j(a)
	\big]
	\kern-2pt
	\big[
	\sigma_k(c)-\sigma_k(b)
	\big]
	-
	\big[
	\sigma_j(c)-\sigma_j(b)
	\big]
	\kern-2pt
	\big[
	\sigma_k(a')-\sigma_k(a)
	\big]
	\kern-3pt
	\Big\}
	h_j
	\wedge
	h_k
= \nonumber \\
	& & 
	\displaystyle
	=
	\sum_{1\le j< k \le n}
	\left\{
	\Big\{
	\big[
	s_{j,k}(a')-s_{j,k}(a)
	\big]
	\wedge
	\big[
	s_{j,k}(c)-s_{j,k}(b)
	\big]
	\Big\}
	\cdot \mathbb{I}_2
	\right\}
	h_j
	\wedge
	h_k \label{eq:balanced bivector}
\end{eqnarray}

where $s_{j,k}=\sigma_j \ell_1 + \sigma_k \ell_2: \Omega \to \mathbb{E}_2$ are $n\choose 2$ smooth transformations. By Theorem~\ref{thm:I} we have that

\[
	\lim_{
	\begin{array}{c}
		\scriptstyle (a,b,c)\to (x,x,x) \\
		\scriptstyle 	a\wedge b +b\wedge c + c\wedge a \ne 0
	\end{array}
	}
	\kern-10pt
	\frac{\Big\{\big[s_{j,k}(a')-s_{j,k}(a)\big]\wedge \big[s_{j,k}(c)-s_{j,k}(b)\big]\Big\}\cdot \mathbb{I}_2}{2\big(a\wedge b +b\wedge c + c\wedge a\big)\cdot \mathbb{I}_2}
	=
	\big(\nabla\sigma_j(x)\wedge \nabla\sigma_k(x)\big)\cdot \mathbb{I}_2\ .
\]

Then, the thesis follows observing that 

\begin{eqnarray}
	& & 
	\displaystyle
	\phantom{=}
	\partial_{\ell_1} s(x) \wedge \partial_{\ell_2} s(x)
		=
	\left(
	\sum_{j=1}^n
	\partial_{\ell_1}\sigma_j(x) h_j
	\right)
	\wedge
	\left(
	\sum_{k=1}^n
	\partial_{\ell_2}\sigma_k(x) h_k
	\right)
	 = \nonumber \\
	& & 
	\displaystyle
	=
	\sum_{1\le j< k \le n}
	\big[
	\partial_{\ell_1}\sigma_j(x)
	\partial_{\ell_2}\sigma_k(x)
	-
	\partial_{\ell_2}\sigma_j(x)
	\partial_{\ell_1}\sigma_k(x)
	\big]
	h_j
	\wedge
	h_k
	= \nonumber \\
	& & 
	\displaystyle
	=
	\sum_{1\le j< k \le n}
	\Big[
	\big(\nabla\sigma_j(x)\wedge \nabla\sigma_k(x)\big)\cdot \mathbb{I}_2
	\Big]
	h_j
	\wedge
	h_k \ . \ \square \label{eq:tangent bivector}
\end{eqnarray}

\section{The area}\label{sec:area}

\begin{theorem}\label{thm:III}
Let $P$ be a compact polygon contained in the open set $\Omega\subseteq \mathbb{E}_2$; let $s:\Omega\to\mathbb{E}_n$ be a smooth surface; let $\Pi$ be a partition of $P$ into a finite family of non-overlapping\footnote{Two sets are {\bf \color{dgreen} non-overlapping}\index{termes}{Non-overlapping sets} if the interiors of those two sets have empty intersection.} nondegenerate oriented triangles $[a_i,b_i,c_i]$ all balanced in $\Omega$; 

\noindent
let~$\displaystyle||\Pi||=\max_{[a_i,b_i,c_i]\in \Pi}\big\{|\ell_{a_i}|,|\ell_{b_i}|, |\ell_{c_i}|, \big\}$; then,

\begin{equation}\label{eq:AlgII}
\lim_{||\Pi||\to 0}
\frac{1}{4}
\sum_{[a_i,b_i,c_i]\in\Pi}
\Big|
	\big[s(a_i')-s(a_i)\big] 
	\wedge
	\big[s(c_i)-s(b_i)\big] 
\Big|
=
\int_P
\big|\partial_{\ell_1} s(x) \wedge \partial_{\ell_2} s(x) \big|dx\ ,
\end{equation}

\noindent
where $a'_i	=
-\left[
a_i 
+ 2 b_i \frac{\ell_{b_i}\cdot \ell_{a_i}}{|\ell_{a_i}|^2} 
+ 2 c_i \frac{\ell_{c_i}\cdot \ell_{a_i}}{|\ell_{a_i}|^2}
\right]$ 
is a balanced mirror vertex for $[a_i,b_i,c_i]$.
\end{theorem}

In particular, if $s$ is injective, the above integral represent the area of $s(P) \subset \mathbb{E}_n$.

\begin{rem*}\label{rem:triangulations}
	In the hypothesis of the foregoing theorem we do not require the partition~$\Pi$ to be a triangulation\footnote{	A {\bf \color{dgreen} triangulation}\index{termes}{Triangulation of a compact polygon} of $P$ is a partition of P into a finite number of non-overlapping triangles such that no vertex of a triangle is an internal point of a side of another.} of $P$ (indeed, a vertex of triangle may be interior to the side of an adjacent triangle). However, if~$\Pi$ is a triangulation of~$P$, then the images~$s(x)$ of its vertices~$x$ are vertices of a polyhedron inscribed on~$s$.
\end{rem*}

\textit{\textbf{Proof of Theorem~\ref{thm:III}.}}

Without loss of generality, we can suppose that all triangles $[a_i,b_i,c_i]\in \Pi$ are equi-oriented with $\mathbb{I}_2=\ell_1\wedge \ell_2$. Then,

\[
\begin{array}{l}
	\displaystyle 
	\left|
		\sum_{[a_i,b_i,c_i]\in \Pi}
		\frac{1}{4} 
			\Big|
				\big[s(a_i')-s(a_i)\big] 
				\wedge
				\big[s(c_i)-s(b_i)\big] 
			\Big|
		-
		\int_P 
\big|\partial_{\ell_1} s(x) \wedge \partial_{\ell_2} s(x) \big|dx
	\right|
	= \\ \\
	\displaystyle 
	=
	\left|
		\sum_{[a_i,b_i,c_i]\in \Pi}
		\left[
		\frac{1}{4} 
			\Big|
				\big[s(a_i')-s(a_i)\big] 
				\wedge
				\big[s(c_i)-s(b_i)\big] 
			\Big|
		-
		\int_{[a_i,b_i,c_i]}
		\big|\partial_{\ell_1} s{\scriptstyle(x)} \wedge \partial_{\ell_2} s{\scriptstyle(x)} \big|dx
		\right]
	\right|
	= \\ \\
	\displaystyle 
	\le
	\sum_{[a_i,b_i,c_i]\in \Pi}
		\left|
		\frac{1}{4} 
			\Big|
				\big[s(a_i')-s(a_i)\big] 
				\wedge
				\big[s(c_i)-s(b_i)\big] 
			\Big|
		-
		\int_{[a_i,b_i,c_i]}
		\big|\partial_{\ell_1} s{\scriptstyle(x)} \wedge \partial_{\ell_2} s{\scriptstyle(x)} \big|dx
	\right|
	= (\#)
\end{array}
\]

Note that the area of each triangle $[a_i,b_i,c_i]$	 ca be written as~$\displaystyle \frac{|u_{a_i}| |\ell_{a_i}|}{2}$, so
	
\[
\begin{array}{l}
	\displaystyle 
	(\#)
	=
	\kern-10pt
	\sum_{[a_i,b_i,c_i]\in \Pi}
		\left|
		\int_{[a_i,b_i,c_i]}
		\left[
		\frac{
			\Big|
				\big[s(a_i')-s(a_i)\big] 
				\wedge
				\big[s(c_i)-s(b_i)\big] 
			\Big|
		}{2 |u_{a_i}| |\ell_{a_i}|}
		-
		\big|\partial_{\ell_1} s{\scriptstyle(x)} \wedge \partial_{\ell_2} s{\scriptstyle(x)} \big|
		\right]\ dx
		\right|
	= \\ \\
	\displaystyle 
	\le
	\kern-10pt
	\sum_{[a_i,b_i,c_i]\in \Pi}
		\int_{[a_i,b_i,c_i]}
		\left|
		\frac{
			\Big|
				\big[s(a_i')-s(a_i)\big] 
				\wedge
				\big[s(c_i)-s(b_i)\big] 
			\Big|
		}{2 |u_{a_i}| |\ell_{a_i}|}
		-
		\big|\partial_{\ell_1} s{\scriptstyle(x)} \wedge \partial_{\ell_2} s{\scriptstyle(x)} \big|
		\right|\ dx \ .
\end{array}
\]

\noindent
Since $\mathbb{G}_{n \choose 2}$ is a Euclidean space, we have that 
$\Big||V|-|W|\Big|\le |V-W|$ for each~$V,W\in \mathbb{G}_{n \choose 2}$, so

\begin{eqnarray*}
		& & 
		\phantom{\le}
		\left|
		\frac{
			\Big|
				\big[s(a_i')-s(a_i)\big] 
				\wedge
				\big[s(c_i)-s(b_i)\big] 
			\Big|
		}{2 |u_{a_i}| |\ell_{a_i}|}
		-
		\big|\partial_{\ell_1} s{\scriptstyle(x)} \wedge \partial_{\ell_2} s{\scriptstyle(x)} \big|
		\right| = 
		\\
		& & 
	\kern-20pt
		\le
		\left|
		\frac{1}{2 |u_{a_i}| |\ell_{a_i}|}
		\Big\{
			\big[s(a_i')-s(a_i)\big] 
			\wedge
			\big[s(c_i)-s(b_i)\big] 
		\Big\}
		-
		\partial_{\ell_1} s{\scriptstyle(x)} \wedge \partial_{\ell_2} s{\scriptstyle(x)}	
		\right| = \\ 
		& &
	\displaystyle
	\kern-20pt
	=
	\kern-2pt
	\left|
	\sum_{1\le j< k \le n}
	\kern-5pt
	\left\{
	\kern-3pt
	\Big\{
		\frac{1}{\scriptstyle 2 |u_{a_i}| |\ell_{a_i}|}
	\big[
	s_{j,k}{\scriptstyle (a'_i)} - s_{j,k}{\scriptstyle (a_i)}
	\big]
	\kern-3pt
	\wedge
	\kern-3pt
	\big[
	s_{j,k}{\scriptstyle (c_i)} - s_{j,k}{\scriptstyle (b_i)}
	\big]
	-
	\nabla\sigma_j{\scriptstyle (x)}
		\kern-3pt
		\wedge 
		\kern-3pt
	\nabla\sigma_k{\scriptstyle (x)}
	\kern-2pt
	\Big\}
	\cdot \mathbb{I}_2
	\kern-2pt
	\right\}
	h_j
	\kern-3pt
	\wedge
	\kern-3pt
	h_k
	\right| \\
	& &
	\displaystyle
	\kern-20pt
	\le
	\kern-2pt
	\sum_{1\le j< k \le n}
	\kern-2pt
	\left|
		\frac{1}{\scriptstyle 2 |u_{a_i}| |\ell_{a_i}|}
	\big[
	s_{j,k}{\scriptstyle (a'_i)} - s_{j,k}{\scriptstyle (a_i)}
	\big]
	\kern-3pt
	\wedge
	\kern-3pt
	\big[
	s_{j,k}{\scriptstyle (c_i)} - s_{j,k}{\scriptstyle (b_i)}
	\big]
	-
	\nabla\sigma_j{\scriptstyle (x)}
		\kern-3pt
		\wedge 
		\kern-3pt
	\nabla\sigma_k{\scriptstyle (x)}
	\right| \ ,
\end{eqnarray*}

\noindent
by equations~(\ref{eq:balanced bivector}) and~(\ref{eq:tangent bivector}). Since each $[a_i,b_i,c_i]$ is equi-oriented with $\mathbb{I}_2$, and we are dealing with bivectors that are also pseudo-scalars, we can write

\begin{eqnarray*}
	& & 
	\left|
		\frac{1}{2 |u_{a_i}| |\ell_{a_i}|}
	\big[
	s_{j,k}(a'_i) - s_{j,k}(a_i)
	\big]
	\kern-3pt
	\wedge
	\kern-3pt
	\big[
	s_{j,k}(c_i) - s_{j,k}(b_i)
	\big]
	-
	\nabla\sigma_j(x)
		\kern-3pt
		\wedge 
		\kern-3pt
	\nabla\sigma_k(x)
	\right|
	= \\
	& & 
	=
		\left|
		\frac{\Big\{\big[s_{j,k}(a'_i)-s_{j,k}(a_i)\big]\wedge \big[s_{j,k}(c_i)-s_{j,k}(b_i)\big]\Big\}\cdot \mathbb{I}_2}{2\big(a\wedge b +b\wedge c + c\wedge a\big)\cdot \mathbb{I}_2}
		-
		\Big(\nabla\sigma_j(x)\wedge \nabla\sigma_k(x)\Big)
		\cdot
		\mathbb{I}_2
		\right| = \\
		& &
		\le
		\left|
			\Big(\nabla\sigma_j(\bar{a}_i)\wedge \nabla\sigma_k(\bar{a}_i)\Big)
			\cdot
			\mathbb{I}_2
			-
			\Big(\nabla\sigma_j(x)\wedge \nabla\sigma_k(x)\Big)
			\cdot
			\mathbb{I}_2
		\right| + \\
	& & 
	\phantom{\le}
	+
	\big|\nabla\sigma_j(\bar{a}_i)\big|
	\left|
		\frac{O\big(|v_{a_i}|^2\big)}{|\ell_{a_i}|} + \frac{O\big(|\ell_{a_i}-v_{a_i}|^2\big)}{|\ell_{a_i}|}
	\right|
		+
	\big|\nabla\sigma_k(\bar{a}_i)\big|\ \Big|O\big(|u_{a_i}|\big)\Big| + \\
	& & 
	\phantom{\le\ }
	+
	\big|\nabla\sigma_k(\bar{a}_i)\big|
	\left|
		\frac{O\big(|v_{a_i}|^2\big)}{|\ell_{a_i}|} + \frac{O\big(|\ell_{a_i}-v_{a_i}|^2\big)}{|\ell_{a_i}|}
	\right|
		+ 
	\big|\nabla\sigma_j(\bar{a}_i)\big| \ \Big|O\big(|u_{a_i}|\big)\Big| + \\ 
	& & 
	\displaystyle
	\phantom{\le \ \ \ }
	+
	O\big(|u_{a_i}|\big)\
	\left|
		\frac{O\big(|v_{a_i}|^2\big)}{|\ell_{a_i}|} + \frac{O\big(|\ell_{a_i}-v_{a_i}|^2\big)}{|\ell_{a_i}|} 
	\right| \ ,
\end{eqnarray*}

by inequality~(\ref{eq:estimate transf}). Then, if we sum over the $n \choose 2$ indexes $j,k$, and if we integrate all over $P$, we obtain quantities that are infinitesimal with respect to $||\Pi||$. $\square$

\section{The graph of a smooth function}

Let $\{\ell_1,\ell_2\}$ be an ordered orthonormal  basis in the Euclidean space $\mathbb{E}_2$; let $\{h_1,h_2,h_3\}$ be an ordered orthonormal  basis in the three-dimensional Euclidean space $\mathbb{E}_3$. We can consider $\mathbb{E}_2 \subset \mathbb{E}_3$ by identifying $h_1=\ell_1$ and $h_2=\ell_2$. So when we have a smooth function $\psi:\Omega \to \mathbb{R}$, we can consider the following smooth surface $s:\Omega \to \mathbb{E}_3$

\begin{center}
$
s(x)=s(\chi\ell_1 + \chi_2\ell_2)
=
\chi\ell_1 + \chi_2\ell_2 + \psi(\chi\ell_1 + \chi_2\ell_2) h_3
=
x+\psi(x)h_3 \ .
$
\end{center}

In this case we have that 

\begin{eqnarray*}
&\phantom{ = } &
\big[s(a')-s(a)\big] \wedge \big[s(c)-s(b)\big] =\\
& = &
(a'-a)\wedge (c-b) - \Big\{\big[\psi(a')-\psi(a)\big](c-b)-\big[\psi(c)-\psi(b)\big](a'-a)\Big\}\wedge h_3 \ ,
\end{eqnarray*}

\begin{eqnarray*}
\partial_{\ell_1} s{\scriptstyle(x)} \wedge \partial_{\ell_2} s{\scriptstyle(x)}
& = & 
\ell_1\wedge\ell_2 + \partial_{\ell_2} \psi(x) \ell_1\wedge h_3 -\partial_{\ell_1}\psi(x) \ell_2\wedge h_3 =
\\
& = & 
\ell_1\wedge\ell_2 - \big(\nabla\psi(x)\big)^* \ ,
\end{eqnarray*}

that is to say, $\displaystyle \nabla\psi(x) = h_3 - \big(\partial_{\ell_1} s{\scriptstyle(x)} \wedge \partial_{\ell_2} s{\scriptstyle(x)}\big)^{\#}= h_3 - \big(\partial_{\ell_1} s{\scriptstyle(x)} \times \partial_{\ell_2} s{\scriptstyle(x)}\big)$,
and

\begin{eqnarray*}
\big|\partial_{\ell_1} s{\scriptstyle(x)} \wedge \partial_{\ell_2} s{\scriptstyle(x)} \big|
& = & 
\sqrt{1 + \big|\nabla\psi(x)\big|^2}\ .
\end{eqnarray*}

%% file: 07-LocalSchwarzParadox.tex
\chapter{The local Schwarz paradox}\label{cha:local Schwarz}

We have seen in Chapter~\ref{cha:smooth curves} (Proposition~\ref{prop: approxim dot c}) that the vector $\dot{c}(\chi)$ is the limit of inscribed mean vectors

\[
\frac{1}{\beta-\alpha}
\big[c(\beta)-c(\alpha)\big] \ ,
\]

\noindent
as $\alpha$ and $\beta$ converge to $\chi$. In this chapter we will verify that the Schwarz Paradox has the following local formulation: inscribed mean bivectors

\begin{center}
$
\displaystyle
\frac{1}{
\left\langle 
a;b;c
\right\rangle
\cdot \mathbb{I}_2
}
\big\langle 
s(a);s(b);s(c)
\big\rangle
$
\end{center}

\noindent
on a smooth surface $s:\Omega\to \mathbb{E}_n$ (such as circular right cylinder) may not converge to the bivector $\partial_{\ell_1}s(x)\wedge \partial_{\ell_2}s(x)$, as $a$, $b$ and $c$ converge to the point $x\in \Omega$. On the contrary, we have seen in Theorem~\ref{thm:II} that the corresponding inscribed balanced mean bivectors 

\begin{center}
$
\begin{array}{l}
	\displaystyle
	\frac{1}{2\left\langle a;b;c\right\rangle\cdot \mathbb{I}_2}
	\big[s(a')-s(a)\big]
	\wedge
	\big[s(c)-s(b)\big]
	= \\ \\
	\displaystyle
	=
	\frac{1}{\big[ \left\langle a;b;c\right\rangle - \left\langle a';b;c\right\rangle\big]\cdot \mathbb{I}_2}
	\Big[
	\big\langle s(a);s(b);s(c)\big\rangle
	-
	\big\langle s(a');s(b);s(c)\big\rangle
 \Big] \ ,
\end{array}
$
\end{center}

\noindent
always converge\footnote{As the non-degenerate triangles $[a,b,c]$ converge to the point $x$.} to the bivector $\partial_{\ell_1}s(x)\wedge \partial_{\ell_2}s(x)$. In this chapter we will verify it on the double sequences of isosceles triangles proposed by Schwarz to prove the fallacy of Serret's definition of area.

\section{The Schwarz triangles}

Let us consider the circular right cylinder of Example~\ref{exa:cylinder}. Such a surface is smooth, and for each $x=\chi\ell_1 +\chi \ell_2 \in \mathbb{E}_2$

\[
\partial_{\ell_1}s(x)\wedge \partial_{\ell_2}s(x)
=
\rho \cos(\chi_1) h_2\wedge h_3 + \rho \sin(\chi_2) h_3\wedge h_1\ . 
\]

Let us consider the double sequences of oriented triangles $[a_{m,n},b_{m,n},c_{m,n}]$ such that

\begin{center}
$
\displaystyle
\hfil
a_{m,n}
=0=x\ ,
\hfil
b_{m,n}
=\frac{\pi}{m}\ell_1 \ + \ \frac{1}{2n}\ell_2\ ,
\hfil
c_{m,n}
=-\frac{\pi}{m}\ell_1 \ + \ \frac{1}{2n}\ell_2\ .
\hfil
$
\end{center}

Then 

\[
x=\lim_{m,n\to \infty} b_{m,n}=\lim_{m,n\to \infty} c_{m,n}=0=a_{m,n}
, \ \ \ \ 
\langle a_{m,n};b_{m,n};c_{m,n} \rangle= \frac{\pi}{mn}\mathbb{I}_2 \ ,
\]
 
\noindent
and

\begin{center}
$
\displaystyle
\big\langle s(a_{m,n});s(b_{m,n});s(c_{m,n}) \big\rangle=
2\rho\sin\frac{\pi}{m}
\left[
\rho\left(1-\cos\frac{\pi}{m}\right)h_1\wedge h_2 + \frac{1}{2n}h_2\wedge h_3
\right]
$
\end{center}

\noindent
so that $
\displaystyle
\frac{1}{
\left\langle 
a_{m,n};b_{m,n};c_{m,n}
\right\rangle
\cdot \mathbb{I}_2
}
\big\langle 
s(a_{m,n});s(b_{m,n});s(c_{m,n})
\big\rangle
$ is asymptotically equivalent to

\begin{center}
$
\displaystyle
2\rho^2n\frac{\pi^2}{m^2}h_1\wedge h_2
+
\rho
h_2\wedge h_3 \ ,
$
\end{center}

\noindent
and then

\begin{itemize}
	\item $\displaystyle
	\kern-5pt
\lim_{m\to\infty}
\frac{1}{
\left\langle 
a_{m,m};b_{m,m};c_{m,m}
\right\rangle
\kern-3pt
\cdot 
\kern-3pt
\mathbb{I}_2
}
\big\langle 
s(a_{m,m});s(b_{m,m});s(c_{m,m})
\big\rangle
\kern-2pt 
=
\kern-2pt 
\rho h_2 \kern-1pt \wedge \kern-1pt h_3
\kern-2pt 
=
\kern-2pt 
\partial_{\ell_1}s{\scriptstyle (0)}\wedge \partial_{\ell_2}s{\scriptstyle (0)}\ ,
$
\item
$\displaystyle
\
\kern-8pt 
\lim_{m\to\infty}
\kern-2pt 
\frac{1}{
\left\langle 
a_{m,m^2};b_{m,m^2};c_{m,m^2}
\right\rangle
\kern-2pt 
\cdot 
\kern-2pt 
\mathbb{I}_2
}
\kern-1pt 
\big\langle 
s(a_{m,m^2});s(b_{m,m^2});s(c_{m,m^2})
\big\rangle
\kern-3pt 
=
\kern-2pt 
2\rho^2 \pi^2 h_1 \kern-2pt \wedge \kern-2pt h_2
+
\rho h_2 \kern-2pt \wedge \kern-2pt h_3 ,
$
\item
and the normalized direction of the mean bivector 

$
\displaystyle
\hfil
\frac{1}{
\left\langle 
a_{m,m^3};b_{m,m^3};c_{m,m^3}
\right\rangle
\cdot \mathbb{I}_2
}
\left\langle 
s(a_{m,m^3});s(b_{m,m^3});s(c_{m,m^3})
\right\rangle
\hfil
$

tends to the planar direction $h_1\wedge h_2$ which is orthogonal to the tangent planar direction $h_2\wedge h_3$.
\end{itemize}

On the contrary, if we consider the mirror vertex $\displaystyle a'_{m,n}=\frac{1}{n}\ell_2$ (that is always balanced), we have that the balanced mean bivector

\begin{center}
$
\displaystyle
	\frac{1}{2\left\langle a;b;c\right\rangle\cdot \mathbb{I}_2}
	\big[s(a')-s(a)\big]
	\wedge
	\big[s(c)-s(b)\big]
	=
	\rho \frac{m}{\pi}
	\sin\frac{\pi}{m} h_2 \wedge h_3\ ,
$
\end{center}

\noindent
converges to the tangent bivector $\rho h_2\wedge h_3$.	

\section{A converging mean bivector not balanced}\label{sec:not balanced bivector}

As we have anticipated in remark~\ref{rem:relaxing hypothesis}

As before, we consider the Schwarz's triangles, but ordered as follows

\begin{center}
$
\displaystyle
\hfil
a_{m,n}
=-\frac{\pi}{m}\ell_1 \ + \ \frac{1}{2n}\ell_2\ ,
\hfil
b_{m,n}
=0=x\ ,
\hfil
c_{m,n}
=\frac{\pi}{m}\ell_1 \ + \ \frac{1}{2n}\ell_2\ .
\hfil
$
\end{center}

As point $d_{(a_{m,n},b_{m,n},c_{m,n})}=d_{m,n}$ we choose $\displaystyle d_{m,n}=\frac{2\pi}{m}\ell_1$ which is not a mirror vertex of $[a_{m,n},b_{m,n},c_{m,n}]$. However, the following relation still holds

\begin{center}
$
\displaystyle
2\left\langle a_{m,n};b_{m,n};c_{m,n}\right\rangle
=
\big\langle a_{m,n};b_{m,n};c_{m,n}\big\rangle - \big\langle d_{m,n};b_{m,n};c_{m,n}\big\rangle
=
\frac{2\pi}{mn} \mathbb{I}_2\ .
$
\end{center}

\begin{eqnarray*}
& & 
\big\langle s(a_{m,n});s(b_{m,n});s(c_{m,n})\big\rangle
 -
\big\langle s(d_{m,n});s(b_{m,n});s(c_{m,n})\big\rangle =\\
& = &
\big[s(d_{m,n})-s(a_{m,n})\big]
	\wedge
\big[s(c_{m,n})-s(b_{m,n})\big]= \\
& = &
\left\{
\rho \left[\cos\left(\frac{2\pi}{m}\right)-\cos\left(\frac{\pi}{m}\right)\right]h_1
+
\rho \left[\sin\left(\frac{2\pi}{m}\right)+\sin\left(\frac{\pi}{m}\right)\right]h_2
-
\frac{1}{2n}h_3
\right\} \wedge \\
& & 
\wedge 
\left\{
\rho \left[\cos\left(\frac{\pi}{m}\right)-1\right]h_1
+
\rho \sin\left(\frac{\pi}{m}\right) h_2
+
\frac{1}{2n}h_3
\right\} = \\
& = & 
\rho \frac{1}{2n}
\left[\sin\left(\frac{2\pi}{m}\right)+2 \sin\left(\frac{\pi}{m}\right)\right] h_2\wedge h_3 
+ 
\rho \frac{1}{2n}
\left[1- \cos\left(\frac{2\pi}{m}\right)\right] h_3\wedge h_1  \ ,
\end{eqnarray*}

that is asymptotically equivalent (as $m,n\to \infty$) to the bivector

\[
\rho
\frac{2\pi}{mn} h_2\wedge h_3
+
\rho
\frac{\pi^2}{m^2 n} h_3\wedge h_1\ .
\]

So we can conclude that

\begin{center}
$
\displaystyle
\lim_{m,n\to \infty}
	\frac{1}{2\left\langle a_{m,n};b_{m,n};c_{m,n}\right\rangle\cdot \mathbb{I}_2}
	\big[s(d_{m,n})-s(a_{m,n})\big]
	\wedge
	\big[s(c_{m,n})-s(b_{m,n})\big]
	=
	\rho h_2 \wedge h_3\ .
$
\end{center}